\documentclass[aos,MSNbibl,dvips]{arximspdf}

\doi{10.1214/14-AOS1262} % kopijuoti is laisko
\volume{42}
\issue{6}
\pubyear{2014}
\firstpage{2557}
\lastpage{2585}
\docsubty{FLA}

\makeatletter
\newcommand{\rrvert}{\vert}
\newcommand{\rrVert}{\Vert}
\newcommand{\llvert}{\vert}
\newcommand{\llVert}{\Vert}
\newcommand{\underset}[3]{\mathop{#3}_{#1}^{#2}}
\newcommand{\eqref}[1]{(\ref{#1})}
\newtheorem{teo}{Theorem}
\newtheorem{lem}{Lemma}
\newtheorem{prop}{Proposition}[section]
\newtheorem{cor}{Corollary}
\newtheorem{lem2}{Lemma}[section]
\newproclaim{rem2}{Remark}[section]
\newcommand{\Cov}{\operatorname{Cov}}
\newcommand{\E}{\mathbb{E}}
\renewcommand{\Re}{\operatorname{Re}}
\renewcommand{\Im}{\operatorname{Im}}
\newcommand{\spa}{\overline{\operatorname{span}}}
\newcommand{\sign}{\operatorname{sign}}

\newcommand{\Sob}{\operatorname{Sob}}
\newcommand{\mE}{\mathcal{E}}
\newcommand{\mG}{\mathcal{G}}
\newcommand{\bbH}{\mathbb{H}}
\newcommand{\bbM}{\mathbb{M}}
\newcommand{\bbL}{\mathbb{L}}
\newcommand{\sym}{\mathrm{sym}}
\newcommand{\bx}{\mathbf{x}}
\newcommand{\bX}{\mathbf{X}}
\newcommand{\by}{\mathbf{y}}
\newcommand{\KL}{\operatorname{KL}}
\makeatother

\begin{document}
\begin{frontmatter}

\title{Asymptotic equivalence for regression under~fractional noise\thanksref{T1}}
\runtitle{Regression under fractional noise}

\begin{aug}
\author[A]{\fnms{Johannes}~\snm{Schmidt-Hieber}\corref{}\ead[label=e1]{schmidthieberaj@math.leidenuniv.nl}}
\runauthor{J. Schmidt-Hieber}
\affiliation{Leiden University}
\address[A]{Mathematical Institute\\
Leiden University\\
Niels Bohrweg 1\\
2333 CA Leiden\\
The Netherlands\\
\printead{e1}}
\end{aug}
\thankstext{T1}{Supported by DFG Fellowship SCHM 2807/1-1 and ERC Grant 320637.}

% HISTORY:
\received{\smonth{12} \syear{2013}}
\revised{\smonth{6} \syear{2014}}

% ABSTRACT
%
\begin{abstract}
Consider estimation of the regression function based on a model with
equidistant design and measurement errors
generated from a fractional Gaussian noise process. In previous
literature, this model has been heuristically linked to an experiment,
where the anti-derivative of the regression function is continuously
observed under additive perturbation by a fractional Brownian motion.
Based on a reformulation of the problem using reproducing kernel
Hilbert spaces, we derive abstract approximation conditions on function
spaces under which asymptotic equivalence between these models can be
established and show that the conditions are satisfied for certain
Sobolev balls exceeding some minimal smoothness. Furthermore, we
construct a sequence space representation and provide necessary
conditions for asymptotic equivalence to hold.
\end{abstract}

% KEYWORDS
\end{frontmatter}

\setcounter{footnote}{1}

%s1 #&#
\section{Introduction}\label{secintro}

Suppose we have observations from the regression model
%
%
%e1 #&#
\begin{equation}
Y_{i,n} = f \biggl(\frac{i}n \biggr) + N^H_{i},
\qquad i=1,\ldots,n, \label{eqmod1}
\end{equation}
where $(N^H_i)_{i\in\mathbb{N}}$ denotes a fractional Gaussian noise
(fGN) process with Hurst index $H\in(0,1)$, that is, a stationary
Gaussian process with autocovariance function $\gamma(k)=\frac{1}2
(\llvert k+1\rrvert ^{2H}-2\llvert k\rrvert ^{2H}+\llvert
k-1\rrvert ^{2H})$. This model can be viewed as a prototype of a
nonparametric regression setting under dependent measurement errors.
Corresponding to $H\leq\frac{1}2$ and $H>\frac{1}2$, the noise process
exhibits short- and long-range dependence, respectively. In the case
$H=\frac{1}2$, fGN is just Gaussian white noise.

Although observing \eqref{eqmod1} is the ``more realistic'' model, one
might be temp\-ted to replace \eqref{eqmod1} by a continuous version
which is more convenient to work with as it avoids discretization
effects. Recall the definition of a fractional Brownian motion (fBM)
with Hurst parameter $H\in(0,1)$ as a Gaussian\vspace*{1pt} process $(B_t^H)_{t\geq
0}$ with covariance function $(s,t)\mapsto\Cov(B_s^H,B_t^H)=\frac{1}2
(\llvert t\rrvert ^{2H}+\llvert s\rrvert ^{2H}-\llvert t-s\rrvert
^{2H})$. In Wang \cite{wan}, Johnstone and Silverman \cite{joh1} and
Johnstone \cite{joh2} it has been argued that
%
%
%e2 #&#
\begin{equation}
Y_t = \int_0^t f(u) \,du
+n^{H-1} B_t^H,\qquad t\in[0,1], B^H
\mbox{ a fBM} \label{eqmod2}
\end{equation}
is a natural candidate for a continuous version of \eqref{eqmod1} for
$H\geq1/2$. By projecting $(Y_t)_{t\in[0,1]}$ onto a suitable basis,
one is further interested in an equivalent sequence space
representation $(\ldots,Z_{-1},Z_0,Z_1,\ldots)$, where the weighted
Fourier coefficients of $f$ under additive white noise are observed,
that is,
%
%
%e3 #&#
\begin{equation}
Z_k= \sigma_k^{-1} \theta_k(f)+n^{H-1}
\varepsilon_k,\qquad k\in\mathbb{Z}, \varepsilon_k
\stackrel{\mathrm{i.i.d.}} {\sim} \mathcal{N}(0,1). \label{eqdefseqmod}
\end{equation}
Here, $\theta_k(f)$ denote the Fourier coefficients and $\sigma_k>0$
are weights. Models of type \eqref{eqdefseqmod} have been extensively
studied in statistical inverse problems literature (cf. Cavalier \cite{cav}).

In this work, we investigate these approximations and its limitations
under Le Cam distance (cf. Appendix~E in the supplementary material \cite{supp}  for a summary of the topic). The
Le Cam distance allows to quantify the maximal error that one might
encounter by changing the experiment. Indeed, it controls the largest
possible difference of risks that could occur under bounded loss
functions. Two experiments are said to be asymptotic equivalent, if the
Le Cam distance converges to zero. Therefore, if we can establish
asymptotic equivalence, then replacing \eqref{eqmod1} by \eqref{eqmod2}
or \eqref{eqdefseqmod} is harmless at least for asymptotic statements
about the regression function $f$.

Our main finding is that for $H\in(\frac{1}4, 1)$ the experiments
generated by model \eqref{eqmod1} and model \eqref{eqmod2} are\vspace*{2pt}
asymptotic equivalent for $\Theta$ a periodic Sobolev space with
smoothness index $\alpha>1/2$, if $H\in[\frac{1}2,1)$ and $\alpha
>(1-H)/(H+1/2)+H-1/2$, if $H\in(\frac{1}4,\frac{1}2]$. Moreover, we
show that for any $H\in(0,1)$ asymptotic equivalence does not hold for
$\alpha=1/2$ and any $\alpha<1-H$, proving that the minimal smoothness
requirement $\alpha>1/2$ for $H\in[\frac{1}2,1)$ is sharp in this
sense. The asymptotic equivalence for $H\in(\frac{1}4, \frac{1}2)$ is
surprising and leads to better estimation rates than the heuristic
continuous approximation presented in \cite{joh1}. The case $H\in
(0,\frac{1}4]$ remains open. Since the noise level in \eqref{eqmod2}
and \eqref{eqdefseqmod} decreases with $H$, discretization effects
become more and more dominant. We conjecture that for small $H$
asymptotic equivalence will not hold. For suitable $\sigma_k, \theta
_k(f)$, equivalence between the experiments generated\vspace*{1pt} by model \eqref
{eqmod2} and model \eqref{eqdefseqmod} can be derived for all $H\in
(0,1)$. We find that $\sigma_k \asymp\llvert k\rrvert ^{1/2-H}$.
Generalization of the latter result is possible if the fBM is replaced
by a Gaussian process with stationary increments.

One of the motivations for our work is to extend the asymptotic
equivalence result for regression under independent Gaussian noise (in
our framework the case $H=1/2$). In Brown and Low \cite{bro}, it was
shown that the experiments generated by the standard regression model
\begin{eqnarray*}
Y_{i,n}&=& f \biggl(\frac{i}n \biggr) +\varepsilon_{i,n},
\qquad(\varepsilon_{i,n})_{i=1,\ldots,n}\stackrel{\mathrm{i.i.d.}} {
\sim} \mathcal{N}(0,1), f\in\Theta
\end{eqnarray*}
and the Gaussian white noise model
%
%
%e4 #&#
\begin{eqnarray}
dY_t &=&f(t)\,dt +n^{-1/2} \,dW_t,\qquad t
\in[0,1], W \mbox{ a Brownian motion, } f\in\Theta,\hspace*{-30pt} \label{eqwhitenmodel}
\end{eqnarray}
are asymptotically equivalent, provided that the parameter space
$\Theta\subset L^2[0,1]$ has the approximation property
%
%
%e5 #&#
\begin{equation}
n \sup_{f\in\Theta} \int_0^1
\bigl(f(u)- \overline f_n(u) \bigr)^2\,du \rightarrow0,
\label{eqBLcond}
\end{equation}
with $\overline f_n:= \sum_{j=1}^n f(j/n)\mathbb
{I}_{((j-1)/n,j/n]}(\cdot)$. A natural choice for $\Theta$ would be the
space of H\"older continuous functions with regularity index larger
$1/2$ and H\"older norm bounded by a finite constant. One of the
consequences of our results is that the smoothness constraints for
general Hurst index, simplify to \eqref{eqBLcond} if $H=1/2$.

\subsection*{Organization of the work}
In Section~\ref{secabstractae}, we use
RKHSs as an abstract tool in order to formulate sufficient
approximation conditions for asymptotic equivalence of the experiments
generated by the discrete and continuous regression models \eqref
{eqmod1} and \eqref{eqmod2}. Although these conditions appear
naturally, they are very difficult to interpret. By using a spectral
characterization of the underlying RKHS, we can reduce the problem to
uniform approximation by step functions in homogeneous Sobolev spaces.
This is described in Section~\ref{sec3}. We further mention
some orthogonality properties that reveal the structure of the
underlying function spaces and explain the key ideas of the proofs. The
main results together with some discussion are stated in Section~\ref
{secmainresults}. In this section, we construct a sequence space
representation and prove equivalence with the continuous regression
experiment. This allows to study the ill-posedness induced by the
dependence of the noise. Thereafter, we study necessary conditions and
outline a general scheme for deriving sequence space representation
given a regression model with stationary noise. This scheme does not
require knowledge of the Karhunen--Loeve expansion. Since the
appearance of \cite{bro}, many other asymptotic equivalence results
have been established for related nonparametric problems and there are
various strategies in order to bound the Le Cam distance. We provide a
brief survey in Section~\ref{secdiscussion} and relate the existing
approaches to our techniques. Proofs are mainly deferred to the
Appendix. Parts of the Appendix as well as a short summary of Le Cam
distance and asymptotic equivalence can be found in the supplement \cite{supp}.

\subsection*{Notation}
If two experiments are equivalent, we write $=$~and~$\simeq$ denotes asymptotic equivalence. The operator $\mathcal{F}$,
defined on $L^1(\mathbb{R})$ or $L^2(\mathbb{R})$ (depending on the
context) is the Fourier transform $(\mathcal{F}f)(\lambda)=\int
e^{-i\lambda x}f(x) \,dx$. For the indicator function on $[0,t)$, we
write $\mathbb{I}_t:=\mathbb{I}_{[0,t)}(\cdot)$ (as a function on
$\mathbb{R}$). The Gamma function is denoted by $\Gamma(\cdot)$. For a
Polish space $\Omega$, let $\mathcal{B}(\Omega)$ be the Borel sets.
Further, $\mathcal{C}[T]$ denotes the space of continuous functions on
$T$ equipped with the uniform norm.

%s2 #&#
\section{Abstract conditions for asymptotic equivalence}\label{secabstractae} %Section 2

The goal of this section is to reduce asymptotic equivalence to
approximation conditions (cf. Theorem~\ref{teomain}). For that, tools
from Gaussian processes and RKHS theory are required which are
introduced in a first step.

The concept of a RKHS can be defined via the Moore--Aronszajn theorem.
It states that for a given index set $T\subseteq\mathbb{R}$ and a
symmetric, positive semi-definite function $K\dvtx T\times T\rightarrow
\mathbb{R}$, there exists a unique Hilbert space $(\bbH, \langle\cdot
,\cdot\rangle_{\bbH})$ with:

\begin{longlist}[(ii)]
\item[(i)] $K(\cdot, t) \in\bbH$, $\forall t\in T$,\vspace*{1pt}
\item[(ii)] $\langle f, K(\cdot,t)\rangle_\bbH= f(t)$, $\forall f\in
\bbH$, $\forall t\in T$.
\end{longlist}
%
% \begin{tabular}{ll}
% (i) & \quad
% $K(\cdot, t) \in\bbH,$ \quad$\forall t\in T$ \\[0.3cm]
% (ii) & \quad$\langle f, K(\cdot,t)\rangle_\bbH= f(t), \quad\forall f
% \end{tabular}

The second condition is called the reproducing property. The Hilbert
space $\bbH$ is called the reproducing kernel Hilbert space (with
reproducing kernel $K$). RKHSs are a strict subset of Hilbert spaces,
as, for instance, $L^2[0,1]$ has no reproducing kernel.

A centered Gaussian process $(X_t)_{t\in T}$ can be associated with a
RKHS $\bbH$ defined via the reproducing kernel $K\dvtx T\times T\rightarrow
\mathbb{R}, K(s,t):=\E[X_sX_t]$. $\bbH$ can also be characterized by
the completion of the function class
%
%
%e6 #&#
\begin{eqnarray}
&& \Biggl\{\phi\dvtx T\rightarrow\mathbb{C} \Big| \phi\dvtx t \mapsto\sum
_{j=1}^M u_j K(s_j,t),
(s_j, u_j) \in T \times\mathbb{C}, j=1,\ldots,M \Biggr
\} \label{eqcomplfctclass}
\end{eqnarray}
with respect to the norm $\llVert \sum_{j=1}^M u_jK(s_j,\cdot)\rrVert
_{\bbH}^2:=\sum_{j,k\leq M} u_j K(s_j,s_k)\overline{u_k}$. For a
Gaussian process $(X_t)_{t\in T}$, there is a Girsanov formula with the
associated RKHS $\bbH$ playing the role of the Cameron--Martin space.

%
%le1 #&#
\begin{lem}[(Example 2.2, Theorem 2.1 and Lemma 3.2 in van der Vaart and van Zanten \cite{vdv2})]
\label{lemcom} %Lemma 1
Let $(X_t)_{t\in T}$ be a Gaussian process with continuous sample
paths on a compact metric space $T$ and $\bbH$ the associated RKHS.
Denote by $P_f$ the probability measure of $t\mapsto f(t)+X_t$ on
$(\mathcal{C}[T], \mathcal{B}(\mathcal{C}[T]))$. If $f\in\bbH$, then
$P_{f}$ and $P_0$ are equivalent measures and the Radon--Nikodym
derivative is given by
\begin{eqnarray*}
\frac{dP_f}{dP_0}&=& \exp \biggl(Uf-\frac{1}2 \llVert f\rrVert
_{\bbH}^2 \biggr),\qquad f\in\bbH,
\end{eqnarray*}
where $U$ denotes the iso-Gaussian process, that is, the centered
Gaussian process $(Uh)_{h\in\mathbb{\bbH}}$ with $UK(t,\cdot):=X_t$
and covariance $\E[(Uh)(Ug)]=\langle h,g\rangle_{\bbH}$.
\end{lem}

Given such a change of measure formula, it is straightforward to
compute the Kullback--Leibler divergence $d_{\KL}(\cdot,\cdot)$ in
terms of the RKHS norm.

%
%le2 #&#
\begin{lem}
\label{lemKLbd}
For $f,g \in\bbH$, and $P_f,P_g$ as in Lemma~\ref{lemcom},
\[
d_{\KL}(P_f,P_g)=\tfrac{1}2 \llVert
f-g\rrVert _{\bbH}^2.
\]
\end{lem}

Throughout the following, the RKHS with reproducing kernel
$(s,t)\mapsto K(s,t)=\E[B_s^HB_t^H]$ will be denoted by $(\bbH,
\llVert \cdot\rrVert _{\bbH})$ (for convenience, the dependence of
$\bbH$ and $K$ on the Hurst index $H$ is omitted). Before the main
result of this section can be stated, we need to introduce the
experiments generated by the models in Section~\ref{secintro}.
\begin{longlist}
\item[\textit{Experiment} $\mE_{1,n}(\Theta)$:]
\textit{Nonparametric regression under fractional noise}. Denote by $\mE_{1,n}(\Theta)=(\mathbb{R}^n,\mathcal
{B}(\mathbb{R}^n),(P_f^n\dvtx  f\in\Theta))$ the experiment with $P_f^n$
the distribution of $\mathbf{Y_n}:=(Y_{1,n},\ldots,Y_{n,n})^t$, where
%
%
%e7 #&#
\begin{equation}
Y_{i,n} = f \biggl(\frac{i}n \biggr) + N^H_{i},
\qquad i=1,\ldots,n\mbox{ and } \bigl(N^H_i
\bigr)_i\mbox{ is a fGN process.} \label{eqmE1nmodel}
\end{equation}
\end{longlist}

\begin{longlist}
\item[\textit{Experiment} $\mE_{2,n}(\Theta)$:]
Let $\mE_{2,n}(\Theta)=(\mathcal
{C}[0,1], \mathcal{B}(\mathcal{C}[0,1]), (Q_f^n\dvtx  f\in\Theta))$ be the
experiment with $Q_f^n$ the distribution of
%
%
%e8 #&#
\begin{equation}
Y_t = \int_0^t f(u) \,du
+n^{H-1} B_t^H,\qquad t\in[0,1], B^H
\mbox{ a fBM}. \label{eqmod2formaldef}
\end{equation}

We write $F_f$ for the anti-derivative of $f$ on $[0,1]$, that is,
$F_f(t)=\int_0^t f(u)\,du$ for all $t\in[0,1]$. The first result relates
asymptotic equivalence to abstract approximation conditions.
\end{longlist}

%
%th1 #&#
\begin{teo}
\label{teomain} %Theorem 1
Let $H\in(0,1)$. Suppose that:
\begin{longlist}[(ii)]
\item[(i)] $(n^{1-2H}\vee1) \sup_{f\in\Theta}\sum_{i=1}^n (n\int_{(i-1)/n}^{i/n} f(u)\,du-f (\frac{i}n ) )^2\rightarrow0$,
\item[(ii)] $n^{1-H} \sup_{f\in\Theta} \inf_{(\alpha_1,\ldots,\alpha
_n)^t\in\mathbb{R}^n} \llVert F_f-\sum_{j=1}^n \alpha_j K(\cdot, \frac
{j}n)\rrVert _{\bbH} \rightarrow0$.
\end{longlist}

Then,
\[
\mE_{1,n}(\Theta)\simeq\mE_{2,n}(\Theta).
\]
\end{teo}

\begin{pf}
The\vspace*{1pt} proof consists of three steps. Proposition~\ref{propE1E5} in the\break
\hyperref[app]{Appendix} states that, under condition~\textup{(i)}, the values $f(\frac{i}n)$
may be replaced by $\widetilde f_{i,n}:=n\int_{(i-1)/n}^{i/n} f(u)\,du$
in model \eqref{eqmE1nmodel}. Instead of $\mE_{1,n}(\Theta)$, we can
therefore consider the experiment $\mE_{4,n}(\Theta)=(\mathbb{R}^n,
\mathcal{B}(\mathbb{R}^n), (P_{4,f}^{n}\dvtx  f\in\Theta))$ with
$P_{4,f}^{n}$ the distribution of
\begin{eqnarray*}
&&\widetilde Y_{i,n}:= \widetilde f_{i,n} +N^H_{i},
\qquad i=1,\ldots,n, f\in\Theta.
\end{eqnarray*}
In order to link experiment $\mE_{4,n}(\Theta)$ to the continuous model
in $\mE_{2,n}(\Theta)$, the crucial point is to construct a path on
$[0,1]$ from the observations $\widetilde Y_{i,n}$, $i=1,\ldots,n$ with
distribution ``close'' to \eqref{eqmod2formaldef}. For this, let
throughout the following $\bx_n=(x_1,\ldots,x_n)$, and $\by
_n=(y_1,\ldots,y_n)$ be vectors and consider the interpolation function
\begin{eqnarray*}
L(t|\bx_n)&:=& \E \bigl[B_t^H |
B_{\ell/n}^H=x_\ell, \ell=1,\ldots,n \bigr],\qquad
t\in[0,1],
\end{eqnarray*}
with $(B_t^H)_{t\geq0}$ a fBM. Let $\mathbf{B}_n^H$ denote the vector
$(B_{\ell/n}^H)_{\ell=1,\ldots,n}$. From the formula for conditional
expectations of multivariate Gaussian random variables, we obtain the
alternative representation
%
%
%e9 #&#
\begin{equation}
L(t| \bx_n) = \biggl(K \biggl(t,\frac{1}n \biggr),K
\biggl(t,\frac{2}n \biggr),\ldots,K (t,1 ) \biggr) \Cov \bigl(
\mathbf{B}_n^H \bigr)^{-1} \bx_n^t
\label{eqLrepresentation}
\end{equation}
and it is easy to verify that
\begin{eqnarray*}
\mbox{\textit{linearity}:}&\qquad& L(\cdot| \bx_n+\by_n)=L(
\cdot| \bx_n)+L(\cdot| \by_n),
\\
\mbox{\textit{interpolation}:}&\qquad& L \biggl(\frac{j}n \bigg| \bx_n
\biggr)= x_j\qquad\mbox{for } j\in\{0,1,\ldots,n\}\mbox{ with }
x_0:=0.
\end{eqnarray*}
The key observations is that if $B^H$ and $\check B^H$ are two
independent fBMs and $R_t^H = \check B_t^H - L(t | (\check B_{\ell
/n}^H)_{\ell=1,\ldots,n})$, then, by comparison of the covariance
structure, the process
\begin{eqnarray*}
&&\bigl(L \bigl(t | \bigl(B_{\ell/n}^H \bigr)_{\ell=1,\ldots,n}
\bigr)+R_t^H \bigr)_{t\geq0}
\end{eqnarray*}
is a fBM as well. Define the vector of partial sums $\mathbf
{S_n\widetilde Y}:=(S_{k}\widetilde Y)_{k=1,\ldots,n}$ with components
$S_{k} \widetilde Y:= \sum_{j=1}^k \widetilde Y_{j,n}$. Recall that
$F_f(t)=\int_0^t f(u)\,du$, let
%
%
%e10 #&#
\begin{equation}
\mathbf{F}_{f,n}:= \biggl(F_f \biggl(\frac{\ell}n
\biggr) \biggr)_{\ell=1,\ldots,n}, \label{eqmathbfFdef}
\end{equation}
and observe that $\mathbf{S_n\widetilde Y}=n \mathbf{F}_{f,n}+ \mathbf
{B}_n^H$, in distribution. For $(R_t)_{t\geq0}$ independent of~$\mathbf
{S_n\widetilde Y}$, we find using the linearity property of $L$,
\begin{eqnarray*}
&&\widetilde Y_t:=n^{H-1} \bigl(L \bigl(t | n^{-H}
\mathbf{S_n\widetilde Y} \bigr)+R_t^H \bigr) =
L(t | \mathbf{F}_{f,n})+ n^{H-1} B_t^H,
\qquad t\in[0,1],
\end{eqnarray*}
for a fBM $B^H$. Consequently, we can construct paths $(\widetilde
Y_t)_{t\in[0,1]}$ by interpolation of $\widetilde Y_{i,n}$, $i=1,\ldots,n$ and adding an uninformative process that match \eqref
{eqmod2formaldef} up to the regression function. On the\vspace*{2pt} contrary, by
the interpolation property of $L$, we can recover $\widetilde Y_{i,n}$,
$i=1,\ldots,n$ from $(\widetilde Y_t)_{t\in[0,1]}$ and, therefore,
\begin{eqnarray*}
\mE_{4,n}(\Theta)&=&\mE_{5,n}(\Theta),
\end{eqnarray*}
where $\mE_{5,n}(\Theta)=(\mathcal{C}[0,1],\mathcal{B}(\mathcal
{C}[0,1]),(Q_{5,f}^n\dvtx  f\in\Theta))$ and $Q_{5,f}^n$ denotes the
distribution of $(\widetilde Y_t)_{t\in[0,1]}$. In Proposition~\ref{propE5E2}, we prove that $\mE_{5,n}(\Theta)\simeq\mE_{2,n}(\Theta)$
under the approximation condition~\textup{(ii)}. This shows that
%
%
%e11 #&#
\begin{equation}
\label{eqoutlineproofofmain} \mE_{1,n}(\Theta) \underset{\mathrm{cond.~(i)}} {
\mathrm{Prop.}~{\fontsize{8.36} {8.36}\selectfont{\ref{propE1E5}}}} {\simeq}
\mE_{4,n}(\Theta) = \mE_{5,n}( \Theta) \underset{
\mathrm{cond.~(ii)}} {\mathrm{Prop.}~{\fontsize{8.36} {8.36}\selectfont{\ref{propE5E2}}}}
{\simeq} \mE_{2,n}(\Theta).
\end{equation}\upqed
\end{pf}

Theorem~\ref{teomain} reduces proving asymptotic equivalence to
verifying the imposed approximation conditions. Whereas \textup{(i)} is of
type \eqref{eqBLcond} and well studied, the second condition requires
that the anti-derivative of $f$ can be approximated by linear
combinations of the kernel functions in the RKHS $\bbH$. In particular,
it implies that $\{F_f\dvtx  f\in\Theta\}\subset\bbH$. In Section~\ref
{secnecessarycond} below, we give a heuristic, why the second condition
appears naturally.

By\vspace*{1pt} Jensen's inequality, property \textup{(i)} in Theorem~\ref{teomain} is
satisfied, provided that $(n^{1-H}\vee n^{1/2})\sup_{f\in\Theta
}\llVert f-\overline f_n\rrVert _{L^2[0,1]}\rightarrow0$ with
$\overline f_n$ being the step function $\sum_{j=1}^n f(\frac
{j}n)\mathbb{I}_{((j-1)/n,j/n]}(\cdot)$. In the case of Brownian
motion, that is $H=1/2$, we can simplify the conditions further. Recall
that in this case $\llVert h\rrVert _{\bbH}=\llVert h'\rrVert
_{L^2[0,1]}$ and $K(s,t)=s\wedge t$ (cf. van der Vaart and van Zanten
\cite{vdv2},\break Section~10). Consequently, both approximation conditions
hold if\break $n^{1/2}\sup_{f\in\Theta}\llVert f-\overline f_n\rrVert
_{L^2[0,1]}\rightarrow0$. Thus,\vspace*{1pt} we reobtain the well-known Brown and
Low condition \eqref{eqBLcond}.

The approximation conditions do not allow for straightforward
construction of function spaces $\Theta$ on which asymptotic
equivalence holds. In view of condition~\textup{(ii)} in Theorem~\ref{teomain}, a
natural class of functions to study in a first step would consists of
all $f$ such that $F_f=K(\cdot,x_0)$ with $x_0\in[0,1]$ fixed, or
equivalently $f\dvtx  t\mapsto\partial_t K(t,x_0)$. In the case $H=\frac
{1}2$, this is just the class of indicator functions $\{\mathbb
{I}_s\dvtx s\in[0,1]\}$ and it is not difficult to see that $\mE
_{2,n}(\Theta)$ is strictly more informative than $\mE_{1,n}(\Theta)$,
implying $\mE_{1,n}(\Theta)\not\simeq\mE_{2,n}(\Theta)$.

Thus, we need to construct $\Theta$ containing smoother functions,
which at the same time can be well approximated by linear combinations
of kernel functions in the sense of condition~\textup{(ii)} of the preceding
theorem. In order to find suitable function spaces, a refined analysis
of the RKHS $\bbH$ is required. This will be the topic of the next section.

%s3 #&#
\section{The RKHS associated to fBM}\label{sec3}

Using RKHS theory, we show in this section that condition~\textup{(ii)} of
Theorem~\ref{teomain} can be rewritten as approximation by step
functions in a homogeneous Sobolev space.

The RKHS of fBM can either be characterized in the time domain via
fractional operators, or in the spectral domain using Fourier calculus.
For our approach, we completely rely on the spectral representation as
it avoids some technical issues. In principle, however, all results
could equally well be described in the time domain. For more on that, cf.
Pipiras and Taqqu \cite{pip}. Set $c_H:=\sin(\pi H)\Gamma(2H+1)$.
Recall that $K(s,t)=\E[B_s^HB_t^H]$, for $s,t\in[0,1]$. Then (cf.
Yaglom \cite{yag} or Samorodnitsky and Taqqu \cite{samo}, equation~(7.2.9)),
%
%
%e12 #&#
\begin{equation}
K(s,t) = \int\mathcal{F}(\mathbb{I}_s) (\lambda)\overline{
\mathcal{F}(\mathbb{I}_t) (\lambda)} \mu(d\lambda) \label{eqKsteq}
\end{equation}
with
\[
\mu(d\lambda)=\frac{c_H}{2\pi} \llvert \lambda\rrvert ^{1-2H}\,d
\lambda.
\]
Given this representation, it is straightforward to describe the
corresponding RKHS as follows (cf. Grenander \cite{gre}, page~97): let
$\bbM$ denote the closed linear span of $\{\mathcal{F}(\mathbb{I}_t)\dvtx
t\in[0,1]\}$ in the weighted $L^2$-space $L^2(\mu)$, then
\begin{eqnarray*}
\bbH&=& \bigl\{F\dvtx \exists F^* \in\bbM,\mbox{ such that } F(t)= \bigl\langle
F^*, \mathcal{F}(\mathbb{I}_t) \bigr\rangle_{L^2(\mu)}, \forall
t\in[0,1] \bigr\},
\end{eqnarray*}
where $\langle g,h \rangle_{L^2(\mu)}:= \int g(\lambda) \overline
{h(\lambda)}\mu(d\lambda)$ denotes the $L^2(\mu)$ inner product. Further,
\[
Q\dvtx \bigl(\bbH,\langle\cdot, \cdot\rangle_{\bbH}\bigr) \rightarrow
\bigl(\bbM,\langle \cdot, \cdot\rangle_{L^2(\mu)}\bigr),\qquad Q(F)=F^*
\]
is an isometric isomorphism and
%
%
%e13 #&#
\begin{equation}
%(13)
\langle g, h\rangle_{\bbH} = \bigl\langle Q(g), Q(h) \bigr
\rangle_{L^2(\mu)}. \label{eqisomspec}
\end{equation}
Let us show the use of this representation of $\bbH$ for the
approximation condition~\textup{(ii)} of Theorem~\ref{teomain}, that is,
%
%
%e14 #&#
\begin{equation}
n^{1-H} \sup_{f\in\Theta} \inf_{(\alpha_1,\ldots,\alpha_n)^t\in\mathbb{R}^n}
\Biggl\llVert F_f-\sum_{j=1}^n
\alpha_j K \biggl(\cdot, \frac{j}n \biggr) \Biggr\rrVert
_{\bbH} \rightarrow0. \label{eqapproxcondiiRESTATE}
\end{equation}
By \eqref{eqKsteq}, we obtain $Q(K(\cdot,s))=\mathcal{F}(\mathbb{I}_s)$
and, therefore,
\[
\Biggl\llVert F_f-\sum_{j=1}^n
\alpha_j K \biggl(\cdot, \frac
{j}n \biggr) \Biggr\rrVert
_{\bbH} = \Biggl\llVert Q(F_f)-\sum
_{j=1}^n \alpha_j \mathcal{F}(
\mathbb{I}_{j/n}) \Biggr\rrVert _{L^2(\mu)}.
\]
It is natural to
consider now functions $f$ for which there exists a $g$ with
$Q(F_f)=\mathcal{F}(g)$. Since $Q(F_f)$ lies in $\bbM$, the closure of
the functions $\{\mathcal{F}(\mathbb{I}_s)\dvtx s\in[0,1]\}$, the support
of $g$ must be contained in $[0,1]$. If for any $f$ such a $g$ exists,
\eqref{eqapproxcondiiRESTATE} simplifies further to
\begin{eqnarray*}
&& n^{1-H} \sup_{f\in\Theta} \inf_{(\beta_1,\ldots,\beta_n)^t\in\mathbb{R}^n}
\Biggl\llVert\mathcal{F} \Biggl(g-\sum_{j=1}^n
\beta_j \mathbb{I}_{ ((j-1)/n,j/n ]} \Biggr) \Biggr\rrVert
_{L^2(\mu)} \rightarrow0.
\end{eqnarray*}
Instead of approximating functions $F_f$ by linear combinations of
kernel functions in $\bbH$, we have reduced the problem to
approximation by step functions in a homogeneous Sobolev space. The
difficulty relies in computing $g$ given a function $f$. To see, how
$f$ and $g$ are linked, observe that by the characterization of $\bbH$
above, $Q(F_f)=\mathcal{F}(g)$, and Parseval's theorem [assuming that
$\llvert \cdot\rrvert ^{1-2H}\mathcal{F}(g)\in L^2(\mathbb{R})$ for
the moment],
\begin{eqnarray*}
F_f(t) &=& \bigl\langle\mathcal{F}(g), \mathcal{F}(
\mathbb{I}_t) \bigr\rangle_{L^2(\mu)}
\\
&=& \frac{c_H}{2\pi}\int_{-\infty}^\infty\llvert
\lambda\rrvert^{1-2H} \mathcal{F}(g) (\lambda) \overline{\mathcal{F}(
\mathbb{I}_t) (\lambda)}\,d\lambda
\\
&=&c_H \int_0^t
\mathcal{F}^{-1} \bigl(\llvert\cdot\rrvert^{1-2H}\mathcal{F}(g)
\bigr) (u)\,du, \qquad t\in[0,1],
\end{eqnarray*}
implying
%
%
%e15 #&#
\begin{equation}
f = c_H\mathcal{F}^{-1} \bigl(\llvert\cdot\rrvert
^{1-2H}\mathcal{F}(g) \bigr)| _{[0,1]}. \label{eqfgrelationship}
\end{equation}
Thus, $f$ and $g$ are connected via a Fourier multiplier restricted to
the interval $[0,1]$. For $H=1/2$, we obtain $f=g$ as a special case.
For other values of $H$, the Fourier multiplier acts like a fractional
derivative/integration operator, in particular it is nonlocal.

One possibility to solve for $g$ is to extend the regression function
$f$ to the real line and then to invert the Fourier multiplier \eqref
{eqfgrelationship}. Recall, however, that $g$ has to be supported on
$[0,1]$ and because of the nonlocality of the operator, this strategy
does not lead to ``valid'' functions $g$.

Another possibility is to interpret \eqref{eqfgrelationship} as a
source condition: we construct function spaces $\Theta$ on which
asymptotic equivalence can be established by considering source spaces,
$\mathcal{S}(\Theta)$ say, of sufficiently smooth functions $g$ first
and then defining $\Theta$ as all functions $f$, for which there exists
a source element $g\in\mathcal{S}(\Theta)$ such that \eqref
{eqfgrelationship} holds. Source conditions are a central topic in the
theory of deterministic inverse problems (cf. Engl et  al.  \cite
{eng} for a general treatment and Tautenhahn and Gorenflo \cite{tau}
for source conditions for inverse problems involving fractional
derivatives). A similar construction is employed in fractional
calculus, by defining the domain of a fractional derivative as the
image of the corresponding fractional integration operator (cf., e.g., see Samko et  al.  \cite{sam}, Section~6.1). Although,
thinking about \eqref{eqfgrelationship} as abstract smoothness
condition itself makes things formally tractable, it has the obvious
drawback, that it does not result in a good description of the function
space $\Theta$. We still cannot decide whether all functions of a given
H\"older or Sobolev space are generated by a source condition or not.

Surprisingly, there are explicit solutions to \eqref{eqfgrelationship},
which satisfy some remarkable orthogonality relations, both in
$L^2[0,1]$ and $L^2(\mu)$. For that some notation is required. Denote
by $J_\nu$ the Bessel function of the first kind with index $\nu>0$. It
is well known that the roots of $J_{\nu}$ are real, countable, nowhere
dense, and also contain zero (cf. Watson \cite{wat}). Throughout the
following, let $\cdots<\omega_{-1}<\omega_0:= 0<\omega_1<\cdots$ be the
ordered (real) roots of the Bessel function $J_{1-H}$ (for convenience,
we omit the dependence on the Hurst index $H$). Define the functions
%
%
%e16 #&#
\begin{equation}
%(16)
g_k\dvtx s\mapsto\mathbb{I}_{(0,1)}(s)
\partial_s \int_0^s
e^{i2\omega_k(s-u)} \bigl(u-u^2 \bigr)^{1/2-H} \,du,\qquad k\in
\mathbb{Z}. \label{eqgkdef}
\end{equation}
As we will show below, $\mathcal{F}^{-1}(\llvert \cdot\rrvert
^{1-2H}\mathcal{F}(g_k)) | _{[0,1]}$ equals (up to a constant factor)
\begin{eqnarray*}
&&f_k\dvtx t\mapsto e^{2i\omega_k t}.
\end{eqnarray*}
This provides us with solutions of \eqref{eqfgrelationship}. It is now
natural to expand functions $f$ as nonharmonic Fourier series $f=\sum_{k=-\infty}^\infty\theta_k e^{2i\omega_k\cdot}$ and to study
asymptotic equivalence with the parameter space $\Theta$ being a
Sobolev ball
%
%
%e17 #&#
\begin{eqnarray}
\label{eqThetaHdef17} \qquad && \Theta_H(\alpha,C)
\nonumber
\\[-2pt]
\\[-16pt]
\nonumber
&&\qquad := \Biggl\{f=\sum_{k=-\infty}^\infty
\theta_k e^{2i\omega_k\cdot}\dvtx \theta_k=\overline{
\theta_{-k}}, \forall k, \sum_{k=-\infty}^\infty
\bigl(1+\llvert k\rrvert \bigr)^{2\alpha}\llvert\theta_k\rrvert
^2\leq C^2 \Biggr\}.
\end{eqnarray}
The constraint $\theta_k=\overline{\theta_{-k}}$ implies that $f$ is
real-valued.

\subsection*{Orthogonality properties of $(f_k)_k$ and $(g_k)_k$}
The\vspace*{1pt} advantage of this approach is that \textit{any} $f\in L^2[0,1]$ can be
expanded in a unique way with respect to $(e^{2i\omega_k\cdot})_k$. We
even have the stronger result.

%
%le3 #&#
\begin{lem}
\label{lemRiesz} %Lemma 3
Given $H\in(0,1)$. Then, $(e^{2i\omega_k\cdot})_k$ is a Riesz basis of
$L^2[0,1]$.
\end{lem}

Recall that a Riesz basis is a ``deformed'' orthonormal basis. For these
bases, Parseval's identity only holds up to constants in the sense of
equivalence of norms. This norm equivalence is usually referred to as
frame inequality or near-orthogonality. For more on the topic, cf.
Young \cite{you}, Section~1.8. The proof of Lemma~\ref{lemRiesz} is
delayed until Appendix~\ref{secproofssecRKHSoffBM}. It relies on a
standard result for nonharmonic Fourier series in combination with
some bounds on the zeros $\omega_k$. Using the previous lemma, the
Sobolev balls $\Theta_H(\alpha,C)$ can be linked to classical Sobolev
spaces for integer $\alpha$; cf. Lemma~\ref{lemSobemb}.

Next, let us prove that $f_k$ and $g_k$ are (up to a constant)
solutions of \eqref{eqfgrelationship} and state the key orthogonality
property of $(g_k)_k$. This part relies essentially on the explicit
orthogonal decomposition of the underlying RKHS $\bbH$ due to
Dzhaparidze and van Zanten \cite{dzh} (cf. also Appendix~\ref
{secproofssecRKHSoffBM}).

%
%th2 #&#
\begin{teo}[(Dzhaparidze and van Zanten \cite{dzh}, Theorem 7.2)]
\label{teodzhvzan} %Theorem 2
Recall that $\cdots<\omega_{-1}<\omega_0:=0<\omega_1<\cdots$ are the
ordered zeros of the Bessel function $J_{1-H}$. For $k\in\mathbb{Z}$, define
%
%
%e18 #&#
\begin{eqnarray}
\phi_k(2\lambda) &=& \sqrt{\frac{\pi}{c_H}} 2^{H-1}
\bigl(1+(\sqrt{2-2H}-1)\delta_{k,0} \bigr) e^{i(\omega_k-\lambda)}
\frac{\lambda^H J_{1-H}(\lambda)}{\lambda-\omega_k}, \label{eqnonnormbasis}
\end{eqnarray}
where $\phi_k(2\omega_k):=\lim_{\lambda\rightarrow2\omega_k}\phi
_k(\lambda)$ and $\delta_{k,0}$ is the Kronecker delta. Then, $\{\phi
_k(\cdot)\dvtx  k\in\mathbb{Z}\}$ is an orthonormal basis (ONB) of $\bbM$
and we have the sampling formula
\begin{eqnarray*}
h&=&\sum_{k=-\infty}^\infty a_k h(2
\omega_k) \phi_k\qquad\mbox{for all } h\in\bbM
\end{eqnarray*}
with convergence in $L^2(\mu)$ and
%
%
%e19 #&#
\begin{eqnarray}
\label{eqakdef} a_k^{-1} &:=& \phi_k(2
\omega_k)
\nonumber
\\[-6pt]
\\[-12pt]
\nonumber
&=& \sqrt{\frac{\pi}{c_H}} \times\cases{ \displaystyle
\sqrt{1-H}2^{2H-3/2} \Gamma(2-H)^{-1}, &\quad for $k=0$,
\vspace*{5pt}
\cr
\displaystyle2^{H-1} \omega_k^HJ_{1-H}'(
\omega_k), &\quad for $k\neq0$.}
\end{eqnarray}
Moreover, for any $k$, $a_{k}=a_{-k}$ and there exists a constant
$\overline c_H$, such that
%
%
%e20 #&#
\begin{equation}
\overline c_H^{-1} \bigl(1+\llvert k\rrvert
\bigr)^{1/2-H}\leq\llvert a_k\rrvert\leq\overline
c_H \bigl(1+\llvert k\rrvert \bigr)^{1/2-H}. \label{eqakbehavior}
%Equation (20)
\end{equation}
\end{teo}

Bessel functions have a power series expansion $J_{1-H}(\lambda)=\sum_{k=0}^\infty\gamma_k\* \lambda^{2k+1-H}$, for suitable coefficients
$\gamma_k$. This\vspace*{1pt} allows to show that $\omega_k=-\omega_{-k}$ for all
integer $k$ and to identify $\lambda^HJ_{1-H}(\lambda)$ for $\lambda<0$
with the real-valued function $\sum_{k=0}^\infty\gamma_k \lambda^{2k+1}$.

The previous theorem is stated in a slightly different form than in
\cite{dzh}; see also the proof in Appendix~\ref{secproofssecRKHSoffBM}.
Let us shortly comment on the sampling formula. Equation (8.544) in
Gradshteyn and Ryzhik \cite{gra} states that $\lambda^HJ_{1-H}(\lambda
)=2^{H-1}\Gamma(2-H)^{-1}\lambda\prod_{k=1}^\infty(1-\frac{\lambda
^2}{\omega_k^2})$. Due to $\omega_k=-\omega_{-k}$, the sampling formula
in Theorem~\ref{teodzhvzan} can thus be rewritten as (infinite)
Lagrange interpolation (cf. also Young \cite{you}, Chapter~4). For
$H=1/2$, we have $J_{1/2}(\lambda)=\sqrt{2/(\pi\lambda)}\sin(\lambda)$
and $\omega_k=k\pi$. In this case, the theorem coincides with a shifted
and scaled version of Shannon's sampling formula.

In the following, we describe the implications of the previous theorem
for our analysis. As an immediate consequence of the sampling formula,
we find that $\langle h,\phi_k\rangle_{L^2(\mu)}=\overline{\langle\phi
_k,h\rangle}_{L^2(\mu)}=a_kh(2\omega_k)$ and
%
%
%e21 #&#
\begin{equation}
\label{eqQspecexplicitderiv} \bigl\langle\phi_k, \mathcal{F}(\mathbb{I}_t)
\bigr\rangle_{L^2(\mu)} =a_k\overline{\mathcal{F}(
\mathbb{I}_t) (2\omega_k)}.
\end{equation}
By Lemma \textup{D.2(ii)} (supplementary material \cite{supp}),
%
%
%e22 #&#
\begin{equation}
%(22)
\mathcal{F}(g_k)=c_H'e^{-i\omega_k}
\phi_k\qquad\mbox{with }c_H':=
\frac{\Gamma(3/2-H)\sqrt{c_H}}{1+(\sqrt{2-2H}-1)\delta_{0,k}} \label{eqcHprimedef}
\end{equation}
and $g_k$ as in \eqref{eqgkdef}. The dependence of $k$ on $c_H'$ is
irrelevant and we can therefore treat it as a constant. Since $\int_0^t
e^{2i\omega_ku}\,du=\overline{\mathcal{F}(\mathbb{I}_t)(2\omega_k)}$, we
have the following chain of equivalences
%
%
%e23 #&#
\begin{eqnarray}
\label{eqQspecexplicitforF} %(23)
f_k(t)=e^{2i\omega_kt} &\quad
\Leftrightarrow\quad& F_{f_k}(t)=\overline{\mathcal{F}(
\mathbb{I}_t) (2 \omega_k)}
\nonumber
\\
&\quad\Leftrightarrow\quad& Q(F_{f_k})= a_k^{-1}
\phi_k
\\
&\quad\Leftrightarrow\quad& Q(F_{f_k})= \mathcal{F} \biggl(
\frac{e^{i\omega_k}
g_k}{a_k c_H'} \biggr).
\nonumber
\end{eqnarray}
This finally shows not only that $f_k$ and $e^{i\omega_k} g_k/(a_k
c_H')$ are solutions to \eqref{eqfgrelationship} but has also two
important further implications for our analysis.

%
%le4 #&#
\begin{lem}
\label{lembiorth} %Lemma 4
The function sequences $(e^{2 i\omega_k \cdot})_k$ and $(a_k e^{i\omega
_k} g_k/c_H')_k$ are bi-orthogonal Riesz bases of $L^2[0,1]$.
\end{lem}

\begin{pf}
By Lemma~\ref{lemRiesz}, $(e^{2i\omega_k\cdot})_k$ is a Riesz basis of
$L^2[0,1]$. From \eqref{eqcHprimedef} and~\eqref{eqnonnormbasis},
$\langle g_k, e^{2i\omega_\ell\cdot}\rangle_{L^2[0,1]}=\mathcal
{F}(g_k)(2\omega_\ell) = c_H'e^{-i\omega_k}\phi_k(2\omega_k)\delta
_{k,\ell}$, with $\delta_{k,\ell}$ the Kronecker delta. Consequently,
$(e^{2i\omega_k\cdot})_k$ and $(a_ke^{i\omega_k}g_k/c_H')_k$ are
biorthogonal implying that $(a_ke^{i\omega_k}g_k/c_H')_k$ is a Riesz
basis of $L^2[0,1]$ as well (cf. Young \cite{you}, page~36).
\end{pf}

Notice that if $f=\sum_k \theta_k e^{2i\omega_k\cdot}$, then in
general, $\theta_k \neq
\langle f, e^{2i\omega_k\cdot}\rangle_{L^2[0,1]}$, since the basis
functions are not orthogonal. Thanks to the previous lemma, the
coefficients~$\theta_k$ can be computed from $f$ via
%
%
%e24 #&#
\begin{equation}
\theta_k= \frac{a_k e^{-i\omega_k} }{c_H'}\langle f, g_k
\rangle_{L^2[0,1]}. \label{eqthetakfromf} %Equation (24)
\end{equation}

Moreover, \eqref{eqQspecexplicitderiv} implies the following explicit
characterization of the RKHS $\bbH$.

%th3 #&#
\begin{teo}
\label{teobbHchar} %Theorem 3
\begin{eqnarray*}
&&\bbH= \Biggl\{F\dvtx F(t)=\sum_{k=-\infty}^\infty
\theta_k \overline{\mathcal{F}(\mathbb{I}_t) (2
\omega_k)}, \sum_{k=-\infty}^\infty
\bigl(1+\llvert k\rrvert \bigr)^{1-2H} \llvert\theta_k\rrvert
^2 <\infty \Biggr\}.
\end{eqnarray*}
\end{teo}

\begin{pf}
Since $(\phi_k)_k$ is an ONB of $\bbM$, $F\in\bbH$ if and only if
$Q(F) =\sum_{k=-\infty}^\infty c_k \phi_k$, with $\sum_{k=-\infty
}^\infty\llvert c_k\rrvert ^2<\infty$. By \eqref
{eqQspecexplicitforF}, this is equivalent to $F(t)=\sum_{k=-\infty
}^\infty\theta_k\overline{\mathcal{F}(\mathbb{I}_t)(2\omega_k)}$ with
$\sum_{k=-\infty}^\infty\llvert a_k\theta_k\rrvert ^2<\infty$. The\vspace*{1pt}
result follows from \eqref{eqakbehavior}.
\end{pf}

Recall the definition of $\Theta_H(\alpha,C)$ in \eqref{eqThetaHdef17}
and set
%
%
%e25 #&#
\begin{equation}
\Theta_H(\alpha):= \bigl\{f\dvtx \exists C=C(f)<\infty\mbox{ with }
f\in\Theta_H(\alpha,C) \bigr\}. \label{eqThetaHdef25}
\end{equation}
From the first equivalence in \eqref{eqQspecexplicitforF}, we obtain

%
%co1 #&#
\begin{cor}
\label{corbbHcharintermsoff} %Corollary 1
$\Theta_H(\frac{1}2 -H)\subset\{f\dvtx  \int_0^\cdot f(u)\,du\in\bbH\}$.
\end{cor}

To conclude this section, notice that we have derived a system of
functions $(f_k,g_k,\phi_k)_k$ with $(f_k)_k$ and $(g_k)_k$ solving
\eqref{eqfgrelationship} and nearly orthogonalizing $\Theta$ and its
source space and $\phi_k$ being an ONB of the underlying RKHS $\bbH$.
The simultaneous (near)-orthogonalization of the spaces is the crucial
tool to verify the second approximation condition of Theorem~\ref{teomain} on Sobolev balls. A slightly simpler characterization of the
RKHS $\bbH$ can be given (cf. Picard \cite{pic}, Theorem 6.12), but it
remains unclear whether it can lead to a comparable simultaneous
diagonalization. For more, see the discussion in Section~\ref{secdiscussion}.

%s4 #&#
\section{Asymptotic equivalence: Main results}\label{secmainresults} %Section 4

%s4.1 #&#
\subsection{Asymptotic equivalence between the experiments \texorpdfstring{$\mathcal{E}_{1,n}(\Theta)$}{$\mathcal{E}_{1,n}(Theta)$} and \texorpdfstring{$\mathcal{E}_{2,n}(\Theta)$}{$\mathcal{E}_{2,n}(Theta)$}}

In this section, we state the theorems establishing asymptotic
equivalence between the experiments generated by the discrete
regression model with fractional measurement noise $Y_{i,n}=f(\frac
{i}n)+N_i^H$, $i=1,\ldots,n$ and its continuous counterpart $Y_t=\int_0^t f(u) \,du +n^{H-1}B_t^H$, $t\in[0,1]$.

Proofs are provided in Appendix~C (supplementary material \cite{supp}).

%
%th4 #&#
\begin{teo}
\label{teoAEexplicitcondlargeH} %Theorem 4
Given $H\in[1/2,1)$. Then, for any $\alpha>1/2$,
\begin{eqnarray*}
&&\mE_{1,n} \bigl(\Theta_H(\alpha,C) \bigr)\simeq
\mE_{2,n} \bigl(\Theta_H(\alpha,C) \bigr).
\end{eqnarray*}
\end{teo}

%
%th5 #&#
\begin{teo}
\label{teoAEexplicitcondsmallH} %Theorem 5
Given $H\in(1/4,1/2)$. If $\Theta_H^{(\sym)}(\alpha,C)
= \{ f\in\Theta_H(\alpha,C)\dvtx  f=-f(1-\cdot)\}$, then, for any $\alpha
>(1-H)/(H+1/2)+H-1/2$,
\begin{eqnarray*}
&&\mE_{1,n} \bigl(\Theta_H^{(\sym)}(\alpha,C) \bigr)
\simeq\mE_{2,n} \bigl(\Theta_H^{(\sym)}(\alpha,C)
\bigr).
\end{eqnarray*}
\end{teo}

In Section~\ref{secnecessarycond}, we show that for any $H\in(0,1)$,
asymptotic equivalence fails to hold if $\alpha=1/2$ or if $\alpha
<1-H$. Therefore, for $H\geq1/2$, the restriction $\alpha>1/2$ is
sharp in this sense. For $H<1/2$, it is more difficult to prove
asymptotic equivalence. If $H\in(1/4,1/2]$, the minimal required
smoothness in the previous result is slightly bigger than the lower
bound $1-H$ but always below $3/4$. In the case $H\downarrow1/4$, the
difference between the upper and lower smoothness assumption becomes
arbitrarily small. For more on the case $H\leq1/4$ and the restriction
to $\Theta_H^{(\sym)}(\alpha,C)$ for $H< 1/2$, see Section~\ref{secnecessarycond}.

In the continuous fractional regression model \eqref{eqmod2}, the noise
term is $n^{H-1}\times$ fBM. Observe that the noise level $n^{H-1}$
corresponds to an i.i.d. (regression) model with $N_n:=n^{2-2H}$
observations. Thus, one can think about $N_n$ as effective sample size
of the problem. If $H<1/2$, we find $N_n\gg n$ and if $H>1/2$, $N_n\ll
n$. The reason for that is the different correlation behavior in the
discrete fractional regression model \eqref{eqmod1}. If $H>1/2$, any
two observations in \eqref{eqmod1} are positively correlated, thus
rewriting this as ``independent'' observations, we obtain $N_n\ll n$. On
the contrary, if $H<1/2$, observations are negatively correlated and
errors cancel out, leading to smaller noise level and, therefore, $N_n
\gg n$.

If short-range dependence is present, that is, $H<1/2$, it has been
argued in Johnson and Silverman \cite{joh1}, that for a specific choice
of $\tau$,
%
%
%e26 #&#
\begin{eqnarray}
\widetilde Y_t&=&\int_0^t f(u)
\,du + \tau n^{-1/2} B_t,\qquad t\in[0,1], B\mbox{ a Brownian
motion} \label{eqcontapproxshortrange}
\end{eqnarray}
is a natural continuous approximation of the discrete model \eqref
{eqmod1}. The advantage is that this does not rely on the fGN and might
hold for any stationary noise process with short-range dependence. For
model \eqref{eqmod1}, however, the asymptotically equivalent continuous
model $Y_t=\int_0^tf(u)\,du+n^{H-1}B_t^H$, has the smaller noise level
$n^{H-1}$, implying that \eqref{eqcontapproxshortrange} leads to a loss
of information.

To conclude the discussion, let us relate the Sobolev ellipsoids $\Theta
_H(\alpha,C)$ to classical Sobolev spaces. From that, we can establish
asymptotic equivalence on a space that depends not on the choice of the
basis. For any positive integer $\beta$, define
\begin{eqnarray*}
&& \Sob_H(\beta,\widetilde C)
\\
&&\qquad := \biggl\{f\in L^2[0,1]\dvtx f^{(\beta-1)}\mbox{ is
absolutely continuous and real-valued},
\\
&&\hspace*{41pt} \llVert f\rrVert_{L^2[0,1]}+ \bigl\llVert f^{(\beta)}
\bigr\rrVert _{L^2[0,1]}\leq\widetilde C,\int_0^1f^{(\ell)}(s)
\bigl(s-s^2 \bigr)^{1/2-H} \,ds=0,
\\
&&\hspace*{273pt} \ell=1,\ldots,\beta \biggr\}.
\end{eqnarray*}

%
%le5 #&#
\begin{lem}
\label{lemSobemb} %Lemma 5
Given $H\in(0,1)$. Then, for any positive integer $\beta$ and
$\widetilde C<\infty$, there exists a finite constant $C$, such that
\[
\Sob_H(\beta,\widetilde C)\subset\Theta_H(\beta,C).
\]
\end{lem}

The proof is delayed until Appendix~C (supplementary material \cite{supp}). For \mbox{$H=\frac{1}2$}, the
constraints  $\int_0^1f^{(\ell)}(s)(s-s^2)^{1/2-H} \,ds=0$, $\ell=1,\ldots,\beta$ simplify to the periodic boundary conditions
$f^{(q)}(0)=f^{(q)}(1)$, $q=0,\ldots,\beta-1$. In this case, Lemma~\ref{lemSobemb} is well known; cf. Tsybakov \cite{tsyb}, Lemma~\textup{A.3}. In the
important case $\beta=1$, the constraint in $\Sob_H(\beta,\widetilde
C)$ is satisfied\vspace*{1pt} whenever $f=f(1-\cdot)$. If we restrict further to
these functions, the definition of $\Sob_H(1,\widetilde C)$ does not
depend on the Hurst index~$H$ anymore. As a consequence of Theorem~\ref{teoAEexplicitcondlargeH} and the embedding, we obtain
the following.

%
%co2 #&#
\begin{cor}
Let $H\in[\frac{1}2,1)$. Then, for any finite constant $\widetilde C$,
\begin{eqnarray*}
&&\mE_{1,n} \bigl(\Sob_H(1,\widetilde C) \bigr)\simeq
\mE_{2,n} \bigl(\Sob_H(1,\widetilde C) \bigr).
\end{eqnarray*}
\end{cor}

%s4.2 #&#
\subsection{Construction of equivalent sequence model}

Let\vspace*{1pt} $\Theta_H(\alpha)$ be as in \eqref{eqThetaHdef25} and write $f= \sum_{k=-\infty}^\infty\theta_k e^{2i\omega_k\cdot}$ for a generic element
of $\Theta_H(\alpha)$. Define the experiment $\mE_{3,n}(\Theta_H(\alpha
))=(\mathbb{R}^{\mathbb{Z}}, \mathcal{B}(\mathbb{R}^{\mathbb
{Z}}),(P_{3,f}^n\dvtx  f\in\Theta_H(\alpha)))$. Here, $P_{3,f}^n$ is the
joint distribution of $(Z_k)_{k\geq0}$ and $(Z_k')_{k\geq1}$ with
%
%
%e27 #&#
\begin{equation}
Z_k= \sigma_k^{-1} \Re(\theta_k)
+n^{H-1} \varepsilon_k\quad\mbox{and}\quad
Z_k'= \sigma_k^{-1} \Im(
\theta_k) +n^{H-1} \varepsilon_k',
\label{eqseqspacereprexact}
\end{equation}
$(\varepsilon_k)_{k\geq0}$ and $(\varepsilon_k')_{k\geq1}$ being two
independent vectors of i.i.d. standard normal random variables. The
scaling factors are $\sigma_k:=a_k/\sqrt{2}$ for $k\geq1$ and $\sigma
_0:=a_0$, with $(a_k)_k$ as defined in \eqref{eqakdef}.

%th6 #&#
\begin{teo}
\label{teoequivbyiso} %Theorem 6
$\mE_{2,n}(\Theta_H(\frac{1}2 -H)) = \mE_{3,n}(\Theta_H(\frac{1}2 -H))$.
\end{teo}

The proof relies completely on RKHS theory and can be found in Appendix~C (supplementary material \cite{supp}). To illustrate the result, let us give an informal derivation here.
First, we may rewrite the continuous fractional regression model~\eqref
{eqmod2} in differential form $dY_t=f(t)\,dt+n^{H-1}\,dB_t^H$. Recall the
definition of $g_k$ in \eqref{eqgkdef} and notice that $g_k=\overline
{g_{-k}}$. Now, let $k\geq1$ and consider the random variables
$Z_k=\int(\overline{e^{i\omega_k}g_k(t)+e^{i\omega
_{-k}}g_{-k}(t)})\,dY_t/(\sqrt{2}c_H')$. Using \eqref{eqthetakfromf},
\begin{eqnarray*}
Z_k&:=& \frac{\sqrt{2}}{a_k} \Re(\theta_k)+
\frac{n^{H-1}}{\sqrt{2}}(\eta_k+\eta_{-k}),\qquad k=1,2,\ldots
\end{eqnarray*}
with $\eta_k:=\int\overline{e^{i\omega_k}g_k(t)}\,dB_t^H /c_H'$. From
Pipiras and Taqqu \cite{pip}, equation (3.4), $E[\int h_1(t)\,dB_t^H
\cdot\int h_2(t)\,dB_t^H]=\langle\mathcal{F}(h_1),\mathcal
{F}(h_2)\rangle_{L^2(\mu)}$ and together with \eqref{eqcHprimedef} and
the fact that $(\phi_k)_k$ is an ONB of $\bbM$,
\begin{eqnarray*}
&&\E[\eta_k\eta_\ell]= \langle\phi_\ell,
\phi_k \rangle_{L^2(\mu)}= \delta_{k,\ell},
\end{eqnarray*}
where $\delta_{k,\ell}$ denotes the Kronecker delta. Hence, $\varepsilon
_k=(\eta_k+\eta_{-k})/\sqrt{2} \sim\mathcal{N}(0,1)$, i.i.d. for
$k=1,2,\ldots$ and $Z_k=\sigma_k^{-1}\Re(\theta_k)+n^{H-1}\varepsilon_k$.
Similarly, we can construct~$Z_0$ and $Z_k'$, $k\geq1$. This shows
informally that the continuous model $\mE_{2,n}(\Theta_H(\frac{1}2-H))$
is not less informative than observing \eqref{eqseqspacereprexact}. The
other direction follows from the completeness of $(g_k)_k$.

As an application of the previous theorem, let us study estimation of
$\theta$ in the model
\begin{eqnarray*}
Y_{i,n}&=&\theta+ N_i^H,\qquad i=1,\ldots,n,
\end{eqnarray*}
that is, model \eqref{eqmod1} with $f=\theta$ constant. To estimate
$\theta$, one could consider the average $\widehat\theta=\frac{1}n
\sum_{i=1}^n Y_{i,n} =\theta+n^{1-H} \xi$, with $\xi$ a standard
normal random variable. By Theorem~\ref{teoAEexplicitcondlargeH} and
Theorem~\ref{teoequivbyiso}, we find, however, that for $H\in[\frac
{1}2, 1)$, the model is asymptotic equivalent to observing $\sigma_0Z_0
= \theta+\sigma_0 n^{H-1} \varepsilon_0$. Recall that $\sigma_0$ depends
on $H$. Clearly, $\sigma_0\leq1$ and $\sigma_0=1$ for $H=1/2$. For
$H>1/2$, numerical evaluation shows that $\sigma_0$ is a bit smaller
than $1$ implying that the estimator $\widehat\theta$ can be slightly
improved. Instead of the sample average, the construction of the
asymptotic equivalence uses a weighted sum over the $Y_{i,n}$'s, where
the weights are chosen proportional to $g_0(i/n)$ (excluding $i=n$)
with $g_0(s)=(s-s^2)^{1/2-H}$ (cf. also the proof of Theorem~\ref{teoAEexplicitcondlargeH}). For $H>1/2$, this gives far more weight to
observations close to the boundaries. The choice of the weighting
function is a consequence of the biorthogonality in Lemma~\ref{lembiorth}.

To conclude the section, let us link the sequence model \eqref
{eqseqspacereprexact} to inverse problems. For that recall that by
\eqref{eqakbehavior}, $\sigma_k \propto\llvert k\rrvert ^{1/2-H}$.
In the case of long-range dependence ($H>1/2$), $\sigma_k\rightarrow0$
and the problem is well-posed. The noise level, however, is $n^{H-1}$
which is of larger order than the classical $n^{-1/2}$. The opposite
happens if $H<1/2$. In this case, the noise level is $o(n^{-1/2})$ but
on the same time we face an inverse problem with degree of
ill-posedness $1/2-H$. In order to illustrate the effects of
ill-posedness and noise level, let us study estimation of $f\in\Theta
_H(\beta,C)$ with smoothness index $\beta>0$ known. Then we might use
the Fourier series estimator $\widehat f=\sum_{\llvert k\rrvert \leq
M_n} \widehat\theta_k e^{2i\omega_k\cdot}$, with $M_n$ some cut-off
frequency. Here, $\widehat\theta_k=\sigma_k (Z_k+iZ_k')$, for $k\geq
0$, with $Z_0':=0$. For negative $k$ set $\widehat\theta_{k}=\overline
{\widehat\theta_{-k}}$. By Lemma~\ref{lemRiesz}, $(e^{2i\omega_k\cdot
})_k$ is a Riesz basis for $L^2[0,1]$ and from the frame inequality
\begin{eqnarray*}
&&\E \bigl[\llVert\widehat f-f\rrVert_{L^2[0,1]}^2 \bigr]
\lesssim\E \biggl[\sum_{\llvert k\rrvert \leq M_n} \llvert
\theta_k-\widehat\theta_k\rrvert^2 \biggr]+
\sum_{\llvert k\rrvert >M_n}\llvert\theta_k\rrvert
^2.
\end{eqnarray*}
Choosing $M_n=O(n^{-(1-H)/(\beta+1-H)})$, the rate becomes $n^{-2\beta
(1-H)/(\beta+1-H)}$ in accordance with Wang \cite{wan} and, for $\beta
=2$, $H\geq1/2$, Hall and Hart \cite{HallHart}. Surprisingly, faster
rates can be obtained if $H$ is small. The ill-posedness is
overcompensated by the gain in the noise level.

%s4.3 #&#
\subsection{Necessary conditions}\label{secnecessarycond}

In this section, we provide necessary minimal smoothness assumptions
for asymptotic equivalence.

The result below shows that asymptotic equivalence cannot hold for the
(smaller) Sobolev space $\Theta_H^{(\sym)}(\alpha,C)$ if $\alpha=1/2$
or $\alpha<1-H$. This shows that $\alpha\geq\frac{1}2 \vee(1-H)$ is
necessary in Theorems~\ref{teoAEexplicitcondlargeH}~and~\ref{teoAEexplicitcondsmallH}.

%
%le6 #&#
\begin{lem} For any $C>0$, if $\alpha=1/2$ or if $\alpha<1-H$, then
\[
\mE_{1,n}\bigl(\Theta_H^{(\sym)}(\alpha, C)\bigr)
\not\simeq\mE_{2,n}\bigl(\Theta _H^{(\sym)}(\alpha, C)
\bigr).
\]
\end{lem}

\begin{pf} We discuss the two cases \textup{(I)} $\alpha=1/2$ and \textup{(II)}
$\alpha<1-H$, separately.
\begin{longlist}[(II)]
\item[(I)] Define $f_{0,n}=cn^{-1/2}\sin(\omega_n(2\cdot-1))$ and
$f_{1,n}= cn^{-1/2}\sin(\omega_{2n}(2\cdot-1))$. Because of $\sin(\omega
_k(2\cdot-1))=(e^{-i\omega_k} e^{2i\omega_k \cdot}-e^{i\omega
_k}e^{-2i\omega_k \cdot})/(2i)$, the constant $c$ can and will be
chosen such that $f_{0,n}, f_{1,n}\in\Theta_H^{(\sym)}(\frac{1}2, C)$
for all $n$. From the equivalent sequence space representation \eqref
{eqseqspacereprexact}, and since by \eqref{eqakbehavior}, $\sigma_k
\asymp\llvert k\rrvert ^{1/2-H}$, we see that $f_{0,n}$ and
$f_{1,n}$ are separable in experiment\break $\mE_{2,n}(\Theta_H^{(\sym)}(\frac
{1}2, C))=\mE_{3,n}(\Theta_H^{(\sym)}(\frac{1}2, C))$ with positive\vspace*{1pt}
probability. Recall that $P_f^n$ denotes the distribution of the
observation vector in experiment\break $\mE_{1,n}(\Theta_H^{(\sym)}(\frac
{1}2,C))$. It is enough to show that
%
%
%e28 #&#
\begin{equation}
d_{\KL} \bigl(P_{f_{0,n}}^n, P_{f_{1,n}}^n
\bigr)\rightarrow0, \label{eqtobesatisfiedforlowerbd}
\end{equation}
since this implies that there exists no test in $\mE_{1,n}(\Theta
_H^{(\sym)}(\frac{1}2, C))$ distinguishing $f_{0,n}$ and $f_{1,n}$
asymptotically with positive probability. Let $x_n=(n+\frac
{1}4(1-2H))\pi$ and notice that $\sin(2x_n\frac{\ell} n)=\sin
(2x_{2n}\frac{\ell} n)$ for all\vspace*{1pt} integer $\ell$. With Lemma~\ref{lemKLbdfGN}, Taylor approximation and Lemma~\textup{D.1} in \cite{supp},
\begin{eqnarray*}
&& d_{\KL} \bigl(P_{f_0}^n, P_{f_{1,n}}^n
\bigr)
\\
&&\qquad \lesssim \bigl(n^{2-2H}\vee n \bigr)\max_{j=1,\ldots,n}
\biggl\llvert f_{0,n} \biggl(\frac{j}n \biggr)-f_{1,n}
\biggl(\frac{j}n \biggr) \biggr\rrvert^2
\\
&&\qquad = \bigl(n^{2-2H}\vee n \bigr)\frac{c^2}{n} \max
_{j=1,\ldots,n} \biggl\llvert\sin \biggl(\omega_n
\frac{2j-n}n \biggr)- \sin \biggl(x_n \frac{2j-n}n \biggr)
\\
&&\hspace*{135pt} {} +\sin \biggl(x_{2n}\frac{2j-n}n \biggr)-
\sin \biggl( \omega_{2n} \frac{2j-n}n \biggr) \biggr
\rrvert^2
\\
&&\qquad \lesssim n^{-2H-1}\vee n^{-2}\rightarrow0.
\end{eqnarray*}
Hence, \eqref{eqtobesatisfiedforlowerbd} holds and this completes the
proof for $\alpha=1/2$.

\item[(II)] Let $\alpha<1-H$ and choose $k_n$ as the smallest integer
larger than $\log n$. Define $f_{0,n}=0$ and $f_{1,n}=
cn^{H-1}k_n^{1/2-H}[\sin(\omega_{k_n}(2\cdot-1))-\sin(\omega
_{k_n+n}(2\cdot-1))]$. If $c=c(C)$ is chosen sufficiently small, then
$f_{0,n}, f_{1,n}\in\Theta_H^{(\sym)}(\alpha, C)$. If $f=f_{1,n}$, the
$k_n$th coefficient is $\theta_{k_n}= cn^{H-1}k_n^{1/2-H}e^{-i\omega
_{k_n}}/(2i)$ and in experiment $\mE_{3,n}(\Theta_H^{(\sym)}(\alpha
,C))$, we observe $Z_{k_n}= \sigma_{k_n}^{-1} \Re(\theta_{k_n}) +
n^{H-1}\varepsilon_{k_n}$ and $Z_{k_n}' = \sigma_{k_n}^{-1} \Im(\theta
_{k_n}) + n^{H-1}\varepsilon_{k_n}'$. Due to $\llvert \sigma
_{k_n}\rrvert \asymp k_n^{1/2-H}$ the functions $f_{0,n}$ and
$f_{1,n}$ can be distinguished in experiment $\mE_{2,n}(\Theta_H^{(\sym
)}(\alpha, C))=\mE_{3,n}(\Theta_H^{(\sym)}(\alpha, C))$ with positive
probability. In contrast, by the same argument as for case~\textup{(I)}, we find
$d_{\KL}(P_{f_{0,n}}^n, P_{f_{1,n}}^n)\rightarrow0$ in $\mE
_{1,n}(\Theta_H^{(\sym)}(\alpha, C))$. This shows that asymptotic
equivalence does not hold.\quad\qed
\end{longlist}\noqed
\end{pf}

The previous lemma shows essentially that the approximation condition~\textup{(i)} in Theorem~\ref{teomain}, which controls the discretization
effects of the regression function, is necessary. Next, we give a
heuristic argument explaining why asymptotic equivalence requires also
an approximation condition in the RKHS, that is, why condition~\textup{(ii)}
in Theorem~\ref{teomain} is necessary. Since under condition~\textup{(i)}, $\mE
_{1,n}(\Theta)\simeq\mE_{5,n}(\Theta)$ [cf. \eqref
{eqoutlineproofofmain}], it is sufficient to study asymptotic\vspace*{1pt}
equivalence between $\mE_{5,n}(\Theta)$ and $\mE_{2,n}(\Theta)$. In $\mE
_{5,n}(\Theta)$, we observe $\widetilde Y_t= L(t| \mathbf
{F}_{f,n})+n^{H-1}B_t^H$, $t\in[0,1]$ with $L(\cdot | \mathbf
{F}_{f,n})$ as in \eqref{eqLrepresentation}. With the change of measure
formula in Lemma~\ref{lemcom}, it is not hard to see that for the
likelihood ratio test $\phi_n=\mathbb{I}_{\{dQ_f^n/dQ_0^n\leq1\}}$,
$Q_f^n \phi_n+Q_0^n(1-\phi_n)\leq2\exp(-\frac{1}8 n^{2-2H}\llVert
F_f\rrVert _{\bbH}^2 )$. Thus, in $\mE_{2,n}(\Theta)$, we can
distinguish with positive probability between $f$ and $0$ if
$n^{1-H}\llVert F_f\rrVert _{\bbH}$ is larger than some constant.
With the same argument, we can distinguish with positive probability
between $f$ and $0$ in experiment $\mE_{5,n}(\Theta)$ provided that
$n^{1-H}\llVert L(\cdot | \mathbf{F}_{f,n})\rrVert _{\bbH}$ is
larger than some constant. For asymptotic equivalence, we should have
therefore that $n^{1-H}\llVert F_f-L(\cdot | \mathbf
{F}_{f,n})\rrVert _{\bbH}$ is small uniformly over $\Theta$ and this
is just a reformulation of condition~\textup{(ii)} (cf. also the proof of
Proposition~\ref{propE5E2}).

\subsubsection*{Necessary conditions for $H<1/2$}
Let us derive a heuristic
indicating that for $H<1/2$, asymptotic equivalence cannot hold on the
unrestricted Sobolev ball $\Theta_H(\alpha,C)$. This motivates the use
of $\Theta_H^{(\sym)}(\alpha,C)$ in Theorem~\ref{teoAEexplicitcondsmallH}. Moreover, we give an argument why asymptotic
equivalence fails for $H\leq1/4$. First, recall that from \eqref
{eqoutlineproofofmain}, the discrete regression experiment $\mE
_{1,n}(\Theta)$ is asymptotically equivalent to $\mE_{5,n}(\Theta)$
under approximation condition~\textup{(i)} of Theorem~\ref{teomain}.
Therefore, $\mE_{1,n}(\Theta)\simeq\mE_{5,n}(\Theta)$, whenever $\Theta
$ is a H\"older ball with index larger $1-H$, for example. To study $\mE
_{1,n}(\Theta)\not\simeq\mE_{2,n}(\Theta)$, it is thus sufficient to
show $\mE_{5,n}(\Theta)\not\simeq\mE_{2,n}(\Theta)$. In $\mE
_{5,n}(\Theta)$, we observe $\widetilde Y_t= L(t| \mathbf
{F}_{f,n})+n^{H-1}B_t^H$, $t\in[0,1]$. Using \eqref
{eqLrepresentation}, $L(\cdot | \mathbf{F}_{f,n})$ is a linear
combination of the functions $K(\cdot, \frac{j}n)=\Cov(B_{\cdot}^H,
B_{j/n}^H)=\frac{1}2 (\llvert \cdot\rrvert ^{2H}+\llvert \frac
{j}n\rrvert ^{2H}-\llvert \cdot-\frac{j}n \rrvert ^{2H})$. Thus,\vspace*{1pt} we
can write $L(\cdot | \mathbf{F}_{f,n})= \sum_{j=1}^n \gamma_{j,n}
K(\cdot, \frac{j}n)$ for suitable weights $(\gamma_{j,n})_j$. In the\vspace*{1pt}
continuous regression experiment $\mE_{2,n}(\Theta)$, $Y_t=\int_0^t
f(u)\,du +n^{H-1} B_t^H$, $t\in[0,1]$ is observed. Informally, we can
consider differentials
\begin{eqnarray*}
\mbox{in } \mE_{5,n}(\Theta){:}&\qquad& d\widetilde Y_t =
\partial_t L(t | \mathbf{F}_{f,n}) \,dt+n^{H-1}
\,dB_t^H,\qquad t\in[0,1],
\\[-2pt]
\mbox{in } \mE_{2,n}(\Theta){:}&\qquad& dY_t = f(t)\,dt +
n^{H-1} \,dB_t^H,\qquad t\in[0,1].
\end{eqnarray*}
The experiments $\mE_{5,n}(\Theta)$ and $\mE_{2,n}(\Theta)$ will be
close if $\partial_t L(t | \mathbf{F}_{f,n})$ well approximates $f(t)$.
Notice, however, that for $H<1/2$, the function $t\mapsto\partial_t
L(t | \mathbf{F}_{f,n})= \sum_{j=1}^n \gamma_{j,n} \partial_tK(t, \frac
{j}n)$ has singularities at $t=\frac{j}n$ for $j=0,1,\ldots,n$ and for
$H\leq1/4$, it is not in $L^2$ anymore. More precisely, $t\mapsto
\partial_tK(t, \frac{j}n)$ has a singularities at $t=0$ and $t=\frac
{j}n$. Since $\partial_t L(t | \mathbf{F}_{f,n})$ and $f(t)$ must be of
the same order, we have $\sum_{j=1}^n \gamma_{j,n}=O(1)$. Typically, $|
\gamma_{j,n}| \lesssim1/n$ for $j=1,\ldots,n$ and this downweights the
singular behavior of $\partial_t L(t | \mathbf{F}_{f,n})$ at $j/n$ for
$j=1,\ldots,n$ but not at $t=0$. Since all summands contribute to the
singularity at $t=0$, we~find\vspace*{-3pt}
\begin{eqnarray*}
&& \partial_t L(t | \mathbf{F}_{f,n}) \sim H \sum
_{j=1}^n \gamma_{j,n} t^{2H-1}
\qquad\mbox{for } t\downarrow0.
\end{eqnarray*}
If $a=a_n\downarrow0$, then $\widetilde Y_{a/n} = \frac{1}2 a^{2H}\sum_{j=1}^n \gamma_{j,n} n^{-2H}(1+o(1))+\eta/(na^H)$ and $Y_{a/n} =
o(n^{-1/2}) + \eta/(na^H)$ with $\eta\sim\mathcal{N}(0,1)$. For
$H<1/2$, $n^{-2H}\gg n^{-1}$ and the paths of $(\widetilde Y_t)_t$ and
$(Y_t)_t$ can be distinguished if $a\downarrow0$ not too fast.
Asymptotic\vspace*{1.5pt} equivalence will thus not hold unless we additionally assume
that $\sum_{j=1}^n \gamma_{j,n}=o(n^{2H-1})$. Via~\eqref{eqLrepresentation}, this can be expressed as a constraint on $f$
indicating why restriction of the Sobolev ball $\Theta_H(\alpha,C)$ in
Theorem~\ref{teoAEexplicitcondsmallH} is necessary for $H<1/2$.

Second, let us give a heuristic argument which shows that asymptotic
equivalence fails to hold for $H\leq1/4$. We compare the decay of
Fourier coefficients for $(\widetilde Y_t)_t$ and $(Y_t)_t$. The
absolute values of the Fourier coefficients $q(a,k):= \int_0^1 e^{2\pi
k t} \sign(t-a)\llvert t-a\rrvert ^{2H-1} \,dt = e^{2\pi i k a}\int_{-a}^{1-a} e^{2\pi i kt} \sign(t)\llvert t\rrvert ^{2H-1} \,dt$ decay
like $\llvert k\rrvert ^{-2H}$ if $a\in[0,1]$. In particular, the
Fourier coefficients depend in a nontrivial way on $a$. Write
$p(a,k)=q(a,k)\llvert k\rrvert ^{2H}$. Then the $k$th Fourier
coefficient of $t\mapsto\partial_t L(t | \mathbf{F}_{f,n})$ is
$\llvert k\rrvert ^{-2H} \sum_{j=1}^n \gamma_{j,n} p(\frac{j}n, k)$,
already assuming that \mbox{$\sum_{j=1}^n \gamma_{j,n}=0$}. The decay is
unaffected by the smoothness of $f$. To compute the Fourier
coefficients of the fBM, we find using \cite{pip}, equation~(3.4) that
$\int_0^1 e^{2\pi i kt} \,dB_t^H \sim\mathcal{N}(0, \llVert e^{2\pi i
k\cdot}\rrVert _{H}^2)$. Since $\llVert e^{2\pi i k\cdot}\rrVert
_{H}^2 = c_H/2\pi\int\llvert \mathcal{F}(e^{2\pi i k\cdot})(\lambda
)\rrvert ^2\* \llvert \lambda\rrvert ^{1-2H} \,d\lambda\asymp\llvert
k\rrvert ^{1-2H}$ (for the last approximation consider a neighborhood
of $\lambda=2\pi k$), roughly,\vspace*{-1pt}
%
%e29 #&#
%e30 #&#
\begin{eqnarray}
\mbox{in } \mE_{5,n}(\Theta){:}&\qquad& \int_0^1
e^{2\pi i k t} \,d\widetilde Y_t \approx k^{-2H}
+n^{H-1} k^{1/2-H} \xi_k,
\qquad k=1,2,\ldots,\nonumber
\\[-2pt]
\mbox{in } \mE_{2,n}(\Theta){:}&\qquad& \int_0^1
e^{2\pi i k t} \,d Y_t \approx\int_0^1
e^{2\pi i k t}f(t)\,dt +n^{H-1} k^{1/2-H} \xi_k,
\nonumber
\\[-4pt]
\eqntext{k=1,2,\ldots,}
\end{eqnarray}
where $\xi_k$ are centered, normally distributed random variables with
variance bounded in $k$. If $k\asymp n$, then in $\mE_{5,n}(\Theta)$,
the Fourier coefficients are in first order $n^{-2H}+n^{-1/2}\xi_k$,
whereas if $f$ is smooth, we observe $o(n^{-1/2})+n^{-1/2}\xi_k$ in $\mE
_{2,n}(\Theta)$. If $H\leq1/4$, we can therefore distinguish $\mE
_{5,n}(\Theta)$ and $\mE_{2,n}(\Theta)$. The only possibility to avoid
this is to add further constraints to $\Theta$ that ensure that $\sum_{j=1}^n \gamma_{j,n} q(\frac{j}n,k)$ is small. Since the argument
above applies to any $k\asymp n$ and $q(\frac{j}n,k)$ depends on $k$,
we exclude more and more subspaces. This indicates $\mE_{5,n}(\Theta)
\not\simeq\mE_{2,n}(\Theta)$ and, therefore, also $\mE_{1,n}(\Theta
)\not\simeq\mE_{2,n}(\Theta)$.

%s4.4 #&#
\subsection{Generalization}

Let us shortly remark on possible extensions of our method. First
notice that Theorem~\ref{teomain} relies on the specific
self-similarity properties of fractional Brownian motion and a
straightforward generalization is only partially possible (cf. Remark~\ref{remeasyext} below). Passing from the continuous model to the
sequence space representation, however, can be stated in a much more
general framework.

Generalizing $\mE_{2,n}(\Theta)$, denote by $\mG_{2,n}(\Theta)=(\mathcal
{C}[0,1],\mathcal{B}(\mathcal{C}[0,1]),(Q_{2,f}^n\dvtx  f\in\Theta))$ the
experiment with $Q_{2,f}^n$ the distribution of
%
%
%e31 #&#
\begin{equation}
Y_t= \int_0^t f(u)\,du +
n^{-\beta} X_t,\qquad t\in[0,1]. \label{eqgenmod}
\end{equation}
Here, $f$ is the regression function, $\beta>0$, and $X:=(X_t)_{t\in
[0,1]}$ is a continuous, centered Gaussian process with stationary
increments. In particular, this contains model \eqref{eqmod2} if $X$ is
a fBM and $\beta=1-H$. The aim of this section is to construct an
equivalent sequence space representation for \eqref{eqgenmod}.

With the Karhunen--Loeve expansion of $X$, this can be done in a
straightforward way. The drawback of this approach is that closed form
formulas for the basis functions are known only for some specific
choices of $X$. Therefore, we propose a different construction leading
again to nonharmonic Fourier series. This approach is based on the
one-to-one correspondence between mass distributions of vibrating
strings and certain measures which was developed in Kre{\u\i}n \cite{kre},
de Branges \cite{bra}, Dym and McKean \cite{dym} and Dzhaparidze et
 al.  \cite{dzh2}.

Let us sketch the construction. Recall that $X$ is a continuous,
centered Gaussian process with stationary increments and let
$K_X(s,t)=\Cov(X_s,X_t)$. We have the representation $K_X(s,t)= \int_{-\infty}^\infty\mathcal{F}(\mathbb{I}_s)(\lambda) \overline{\mathcal
{F}(\mathbb{I}_t)(\lambda)}\,d\mu_X(\lambda)$, where $\mu_X$ is a
symmetric Borel measure on $\mathbb{R}$ satisfying $\int_{-\infty
}^\infty(1+\lambda)^{-2}\,d \mu_X(\lambda)<\infty$ (cf. Doob \cite{doo},
Section~XI.11). If $\bbM_X=\spa\{\mathcal{F}(I_t)\dvtx t\in[0,1]\}\subset
L^2(\mu_X)$, then the RKHS $\bbH_X$ associated to the Gaussian process
$X$ is given by
\begin{eqnarray*}
\bbH_X&=& \bigl\{F\dvtx \exists F^* \in\bbM_X, \mbox{
such that } F(t)= \bigl\langle F^*, \mathcal{F}(\mathbb{I}_t) \bigr
\rangle_{L^2(\mu_X)}, \forall t\in[0,1] \bigr\}
\end{eqnarray*}
and we have the isometric isomorphism
\begin{eqnarray*}
Q_X\dvtx \bigl(\bbH_X, \langle\cdot,\cdot
\rangle_{\bbH_X} \bigr)&\rightarrow& \bigl(\bbM_{\mu_X},\langle
\cdot,\cdot\rangle_{L^2(\mu_X)} \bigr),
\\
F&\mapsto& F^*.
\end{eqnarray*}
In order to extend Theorem~\ref{teodzhvzan}, the crucial observation is
that there is a one-to-one correspondence between measures $\mu_X$ with
$\int_{-\infty}^\infty(1+\lambda)^{-2}\,d \mu_X(\lambda)<\infty$ and
mass distribution functions, say $m$, of a vibrating string.
Computation of $m$ is quite involved and thus omitted here. For a
detailed explanation, see \cite{dzh2}. If $m$ is continuously
differentiable and strictly positive, we have the following theorem.

%
%th7 #&#
\begin{teo}
\label{teogendzh}
Denote by $(\lambda, \omega) \mapsto S(\lambda,\omega)$ the
reproducing kernel of $\bbM_X$. There exist real numbers $\cdots<\nu
_{-1}<\nu_0=0<\nu_1<\cdots$ such that $(S(\nu_k, \cdot)/\break \sqrt{S(\nu
_k,\nu_k)})_k$ is an ONB of $\bbM_X$ and
\begin{eqnarray*}
h &=& \sum_{k=-\infty}^\infty h(\nu_k)
\frac{S(\nu_k, \cdot)}{S(\nu_k,\nu_k)}\qquad\mbox{for all } h\in\bbM_X
\end{eqnarray*}
with convergence of the sum in $L^2(\mu_X)$. Furthermore, for any $k\in
\mathbb{Z}$, $\nu_k=-\nu_{-k}$, $S(\nu_k,\nu_k)=S(\nu_{-k},\nu_{-k})$
and $\llvert \nu_k\rrvert =2\llvert k\rrvert \pi(1+o(1))$, for
$\llvert k\rrvert \rightarrow\infty$.
\end{teo}

\begin{pf}
This follows largely from Dzhaparidze et  al.  \cite{dzh2},
Theorem 3.5 and Zareba \cite{zar}, Lemma 2.8.7. It remains to show that
$\nu_k=-\nu_{-k}$ and $S(\nu_k,\nu_k)=S(\nu_{-k},\nu_{-k})$. Notice
that by the reproducing property $S(\nu_k,\nu_k)=\llVert S(\nu
_k,\lambda)\rrVert _{L^2(\mu_X)}^2\geq0$. To see that $\nu_k=-\nu
_{-k}$, observe that from \cite{dzh2}, equation (2.2), $A(x,\lambda
)=A(x,-\lambda)$. Therefore, $B(x,\lambda)=-\frac{1}{\lambda
}A^+(x,\lambda)$ (cf. \cite{dzh2}, Section~2.3) satisfies $B(x,\lambda
)=-B(x,-\lambda)$. The numbers $\cdots<\nu_{-1}<\nu_0=0<\nu_1<\cdots$
are the zeros of $B(x,\cdot)$ for a specific value of $x$ (cf. \cite
{dzh2}, Theorem 2.10, equation (3.2) and Theorem 3.5). Hence, $\nu
_{-k}=-\nu_k$. Using \cite{dzh2}, equation (2.10),
\begin{eqnarray*}
K_T(\nu_k,\lambda) &=& \frac{A(x(T),\nu_k)B(x(T),\lambda)}{\pi(\lambda
-\nu_k)}
\\
&=& \frac{A(x(T),\nu_{-k})B(x(T),-\lambda)}{\pi(-\lambda-\nu
_{-k})}
\\
&=&K_T(\nu_{-k},\lambda).
\end{eqnarray*}
Together with \cite{dzh2}, equation (3.1) this shows $S(\nu_k,\nu
_k)=S(\nu_{-k},\nu_{-k})$. The proof is complete.
\end{pf}

Write $\psi_{k,X}=S(\nu_k, \cdot)/\sqrt{S(\nu_k,\nu_k)}$ and $\rho
_k=S(\nu_k,\nu_k)^{-1/2}$ and notice that by the preceding theorem,
$\rho_k=\rho_{-k}$. The sampling formula reads then $h=\sum_{k=-\infty
}^\infty\rho_k h(\nu_k)\psi_{k,X}$ and $\langle h,\psi_{k,X}\rangle
_{L^2(\mu_X)}=\rho_k h(\nu_k)$. Generalizing\vspace*{1pt} \eqref
{eqQspecexplicitforF},\break we find that if $F(t)=\sum_{k=-\infty}^\infty
\theta_k\mathcal{F}(\mathbb{I}_t)(\nu_k)$, then $F^*=Q_X(F)=\break \sum_{k=-\infty}^\infty\theta_k\rho_k^{-1}\psi_{k,X}$. Analogously to
Theorem~\ref{teobbHchar} and Corollary~\ref{corbbHcharintermsoff}, we
have the characterization
\begin{eqnarray*}
\bbH_X&=& \Biggl\{F\dvtx F(t)=\sum_{k=-\infty}^\infty
\theta_k\overline{\mathcal{F}(\mathbb{I}_t) (
\nu_k)}, \sum_{k=-\infty}^\infty
\rho_k^2\llvert\theta_k\rrvert ^2<
\infty \Biggr\}
\end{eqnarray*}
and that $\Theta_X:=\{f\dvtx  f=\sum_{k=-\infty}^\infty\theta_ke^{i\nu
_k\cdot}, \theta_k=\overline\theta_{-k} \sum_{k=-\infty}^\infty\rho
_k^2\llvert \theta_k\rrvert ^2<\infty\}$ is a subset of $\{f\dvtx  \int_0^{\cdot} f(u)\,du\in\bbH_X\}$. Due to $\nu_k=-\nu_{-k}$, a function is
real-valued within this class iff $\theta_k=\overline{\theta_{-k}}$ for
all $k$.

Generalizing $\mE_{2,n}$, define the experiment $\mG_{3,n}(\Theta
_X)=(\mathbb{R}^{\mathbb{Z}}, \mathcal{B}(\mathbb{R}^{\mathbb
{Z}}),(Q_{3,f}^n\dvtx  f\in\Theta_X))$. Here, $Q_{3,f}^n$ is the joint
distribution of $(Z_k)_{k\geq0}$ and $(Z_k')_{k\geq1}$ with
\begin{eqnarray*}
&& Z_k= \sigma_k^{-1} \Re(\theta_k)
+n^{H-1} \varepsilon_k\quad\mbox{and}\quad
Z_k'= \sigma_k^{-1} \Im(
\theta_k) +n^{H-1} \varepsilon_k',
\end{eqnarray*}
$(\varepsilon_k)_{k\geq0}$ and $(\varepsilon_k')_{k\geq1}$ are two
independent vectors of Gaussian noise. The scaling factors are $\sigma
_k:=\rho_k/\sqrt{2}$ for $k\geq1$ and $\sigma_0:=\rho_0$.

%
%th8 #&#
\begin{teo} $\mG_{2,n}(\Theta_X)=\mG_{3,n}(\Theta_X)$.
\end{teo}

The proof is the same as for Theorem~\ref{teoequivbyiso}.

The advantage of this approach is that by following the program
outlined in \cite{dzh2} or \cite{zar}, closed form expressions for
$\sigma_k$ can be derived even if the Karhunen--Loeve decomposition is
unknown. The difference is that $f$ is not expanded in an ONB but again
as a nonharmonic Fourier series with respect to $(e^{i\nu_k\cdot})_k$.
By Theorem~\ref{teogendzh}, $\llvert \nu_k\rrvert =2\llvert
k\rrvert \pi(1+o(1))$. Therefore,\vspace*{1pt} the functions $(e^{i\nu_k\cdot})_k$
are ``close'' to the harmonic basis $(e^{2\pi i k\cdot})_k$.

%s5 #&#
\section{Discussion}\label{secdiscussion}

In this section, we give a short summary of related work on regression
under dependent noise and asymptotic equivalence.

Optimal rates of convergence for regression under long-range dependent
noise were first considered by Hall and Hart \cite{HallHart} using
kernel estimators.

Inspired by the asymptotic equivalence result of Brown and Low \cite
{bro}, Wang \cite{wan} makes the link between discrete regression under
dependent noise and experiment $\mE_{2,n}(\Theta)$, in which the path
of the integral of $f$ is observed plus a scaled fBM. Passing from the
discrete to the continuous model is done by adding uninformative
Brownian bridges. It is argued that this will lead to good
approximations of the continuous path. From an asymptotic equivalence
perspective this interpolation scheme leads, however, to dependencies
in the errors which are difficult to control. To prove Theorem~\ref{teomain}, we used instead the interpolation function
%
%
%e32 #&#
\begin{equation}
t\mapsto L(t| \bx_n)= \E \bigl[B_t^H |
B_{\ell/n}^H=x_\ell, \ell=1,\ldots,n \bigr],
\label{eqLdiscussed}
\end{equation}
which has the advantage that the interpolated discrete observations
have the exact error distribution of the continuous model resulting in
the equivalence $\mE_{4,n}(\Theta)=\mE_{5,n}(\Theta)$ in the proof of
Theorem~\ref{teomain}. The use of the interpolation function \eqref
{eqLdiscussed} for asymptotic equivalence appears implicitly already in
the proof of Theorem 2.2 in Rei{\ss} \cite{rei2}.

Approximation\vspace*{1pt} of discrete regression under dependent errors by\break $\int_0^\cdot f(u) \,du+n^{H-1} B_\cdot^H$ and sequence model representations
were further studied in Johnstone and Silverman \cite{joh1} and more
detailed in Johnstone \cite{joh2}.

Donoho \cite{Donoho1995101} investigates the wavelet-vaguelette
decomposition for inverse problems. Since this is very close to the
simultaneous orthogonalization presented in Section~\ref{sec3},
the connection is discussed in more detail here. Let $f_k=e^{2i\omega
_k\cdot}$ and $\check g_k = a_ke^{i\omega_k} g_k/c_H'$. By Lemma~\ref{lembiorth},
\begin{eqnarray*}
&& (f_k)_k\quad\mbox{and}\quad(\check g_k)_k
\qquad\mbox{are biorthogonal bases of } L^2[0,1].
\end{eqnarray*}
Next, define the operator $Sh = \sqrt{c_H/2\pi} \llvert \cdot\rrvert
^{1/2-H} \mathcal{F}(h)$ and notice that \eqref{eqfgrelationship} can
be rewritten as $f=S^*Sg$ with $S^*$ the adjoint operator. By Theorem
\ref{teodzhvzan}, the functions $\lambda\mapsto\psi_k(\lambda) =\sqrt
{c_H/2\pi} \llvert \lambda\rrvert ^{1/2-H}\phi_k(\lambda)$ are
orthonormal with respect to $L^2(\mathbb{R})$. Using \eqref
{eqcHprimedef} and $f_k=S^*S (a_k^{-2}\check g_k)$, we have the
quasi-singular value decomposition
\begin{eqnarray*}
S \check g_k &=& a_k\psi_k \quad\mbox{and}
\quad S^*\psi_k = a_k f_k\qquad\mbox{for all }
k\in\mathbb{Z}.
\end{eqnarray*}
This should be compared to the wavelet-vaguelette decomposition in \cite
{Donoho1995101}, Section~1.5 which proposes, within a general
framework, to start with a wavelet decomposition $(\psi_{j,k})_{j,k}$,
replacing the orthonormal functions $(\psi_k)_k$ above. This allows to
consider more general spaces than Sobolev balls but cannot be applied
here as we need to work in the RKHS $\bbH$. The restriction to $\bbH$
implies that the underlying functions $\phi_k$ (or $\phi_{j,k}$ in a
multi-resolution context) have to be functions in $\bbM$, that is, the
closed linear span of $\{\mathcal{F}(\mathbb{I}_t)\dvtx  t\in[0,1]\}$ in
$L^2(\mu)$. Because fBM on $[0,1]$ is considered, the index $t$ has to
be in $[0,1]$ and $\bbM$ is very difficult to characterize. In
particular, it is strictly smaller than $L^2(\mu)$. This shows that
finding an ONB of $\bbM$ (cf. Theorem~\ref{teodzhvzan}) is highly
nontrivial and it remains unclear in which sense $\bbM$ could admit a
multi-resolution decomposition.

Besides that, Brown and Low \cite{bro} and Nussbaum \cite{nus}
established nonparametric asymptotic equivalence as own research
field. Since then, there has been considerable progress in this area.
Asymptotic equivalence for regression models was further generalized to
random design in Brown et  al.  \cite{brown2002}, non-Gaussian
errors in Grama and Nussbaum \cite{grama} and higher-dimensional
settings in Carter \cite{carter2006} and Rei{\ss}~\cite{rei}. Rohde \cite
{roh2} considers periodic Sobolev classes, improving on condition~\eqref{eqBLcond} in this case. Carter \cite{car2009} establishes asymptotic
equivalence for regression under dependent errors. The result, however,
is derived under the strong assumption that the noise process is
completely decorrelated by a wavelet decomposition. Multiscale
representations that nearly whiten fBM are known (cf. Meyer et  al.  \cite{MeyerSellanTaqqu}, Section~7), but it is unclear whether fBM
admits an exact wavelet decomposition. One possibility to extend the
result to regression under fractional noise is to give up on
orthogonality and to deal with nearly orthogonal wavelet decompositions
instead. This, however, causes various new issues that are very
delicate and technical. One might view the methods developed in Golubev
et  al.  \cite{gol} and Rei{\ss}~\cite{rei2} as first steps toward
such a theory, as both deal with similar problems, however in very
specific settings.

%sA #&#
\begin{appendix}\label{app}
%sB #&#
\section{Proofs for Section~\texorpdfstring{\lowercase{\protect\ref{secabstractae}}}{2}}\vspace*{-6pt}\label{secappendixproofsforsecabstractae}

\begin{pf*}{Proof of Lemma~\ref{lemKLbd}}
Write $\nu_f= dP_f/dP_0$ and $\nu_g=dP_g/dP_0$. Moreover, denote by
$E_0[\cdot]$ expectation with respect to $P_0$. We have $E_0[\nu_f]=1$
and by Lemma~\ref{lemcom},
\begin{eqnarray*}
d_{\KL}(P_f,P_g) &=& E_0 \biggl[
\log \biggl(\frac{\nu_f}{\nu_g} \biggr) \nu_f \biggr] =
E_0 \bigl[ (Uf-Ug) \nu_f \bigr]-\frac{1}2
\llVert f\rrVert_{\bbH}^2+ \frac{1}2 \llVert g\rrVert
_{\bbH}^2. % \label{eqKLexp}
\end{eqnarray*}
Note that $g=g_1+g_2$ with $g_1:=\langle g,f\rangle_{\bbH} \llVert
f\rrVert _{\bbH}^{-2}f$ and $g_2:=g-g_1$. Clearly, $\Cov(Uf,Ug_2)=0$
and since $Uf, Ug_2$ are Gaussian, $\nu_f$ and $Ug_2$ are independent.
For a centered normal random variable $\xi$ with variance $\sigma^2$,
\[
\E\bigl[\xi\exp(\xi)\bigr]=\partial_t \E\bigl[\exp(t\xi)\bigr]|
_{t=1}= \partial_t \exp \biggl(\frac{\sigma^2}2
t^2 \biggr) \bigg| _{t=1}=\sigma^2\exp \biggl(
\frac{\sigma^2}2 \biggr)
\]
and hence,
\begin{eqnarray*}
&&E_0 \bigl[ (Uf-Ug) \nu_f \bigr]= \biggl(1-
\frac{\langle g,f\rangle_{\bbH}}{\llVert f\rrVert _{\bbH}^2} \biggr) E_0 \bigl[ (Uf) \nu_f \bigr]=
\llVert f\rrVert_{\bbH}^2-\langle g,f\rangle_{\bbH}.
\end{eqnarray*}
Plugging this into the formula for $d_{\KL}(P_f,P_g)$, the result follows.
\end{pf*}

For a similar result, cf. Gloter and Hoffmann \cite{glo3}, Lemma 8.

%sB.1 #&#
\subsection{Completion of Theorem \protect\texorpdfstring{\ref{teomain}}{1}}

The remaining parts for the proof of Theorem~\ref{teomain} follows from Propositions~\ref{propE1E5}~and~\ref{propE5E2} below.

If $\bX=(X_1,\ldots,X_n)$ is a stationary process with spectral density $f$, it is well known that the
eigenvalues of the Toeplitz matrix $\Cov(\bX)$ lie between the minimum and maximum of $f$ on $[0,\pi]$.
These bounds become trivial if $f(\lambda)$ converges to $0$ or $\infty$ for $\lambda\downarrow 0$. The
first lemma gives a sharper lower bound for the smallest eigenvalue of a Toeplitz matrix which is of independent interest.

%leB.1 #&#
\begin{lem2}
\label{lemlowerbdstatprocesses}
 Let $\bX=(X_1,\ldots,X_n)$ be a stationary process with spectral density $f$ and denote by $\lambda_i(\cdot)$  the $i$th eigenvalue. Then
\begin{eqnarray*}
&&\lambda_n \bigl(\Cov(\bX) \bigr)\geq \biggl(1-\frac{1}{\pi}
\biggr) \inf_{\lambda \in  [1/n, \pi ]} f(\lambda).
\end{eqnarray*}
\end{lem2}

\begin{pf}
For any vector $v=(v_1,\ldots,v_n)$,
\begin{eqnarray*}
&& v^t \Cov(\bX) v = \frac{1}{2\pi} \int_{-\pi}^{\pi}
\Biggl\llvert \sum_{k=1}^n
v_k e^{ik\lambda} \Biggr\rrvert ^2 f(\lambda) \,d
\lambda = \frac{1}{\pi} \int_0^{\pi} \Biggl
\llvert \sum_{k=1}^n v_k
e^{ik\lambda} \Biggr\rrvert ^2 f(\lambda) \,d\lambda.
\end{eqnarray*}
In particular, $\llVert  v\rrVert  ^2=v^tv= \frac{1}{\pi} \int_0^{\pi}  \llvert   \sum_{k=1}^n v_k e^{ik\lambda}  \rrvert  ^2 \,d\lambda$. The estimate,
\begin{eqnarray*}
&& v^t \Cov(\bX) v \geq \Biggl(\llVert v\rrVert ^2-
\frac{1}{\pi}\int_0^{1/n} \Biggl\llvert \sum
_{k=1}^n v_k e^{ik\lambda}
\Biggr\rrvert ^2 \,d\lambda \Biggr) \inf_{\lambda \in  [1/n, \pi ]} f(
\lambda)
\end{eqnarray*}
together with $ \llvert   \sum_{k=1}^n v_k e^{ik\lambda}  \rrvert  ^2\leq  (\sum_{k=1}^n \llvert  v_k\rrvert   )^2\leq n\llVert  v\rrVert  ^2$ yields the result.
\end{pf}
Along the line of the proof, one can also show that $\sup_{\lambda\in [1/n, \pi)}f(\lambda) + \frac{n}{\pi} \int_0^{1/n} f(\lambda ) \,d\lambda$ is an upper bound of the eigenvalues.

%leB.2 #&#
\begin{lem2}
\label{lemKLbdfGN}
 For a vector $v\in \mathbb{R}^n$, let $P_v$ denote the distribution of $(Y_{1,n},\ldots, Y_{n,n})$ with $Y_{i,n}=v_i+N_i^H$, $i=1,\ldots,n$ and $(N^H_i)_i$ fGN. Then there exists a constant $c=c(H)$, such that
\begin{eqnarray*}
&& d_{\KL}(P_v,P_w)\leq c \bigl(n^{1-2H}
\vee 1 \bigr) (v-w)^t(v-w).
\end{eqnarray*}
\end{lem2}

\begin{pf}
Denote the spectral density of fractional Gaussian noise with Hurst index $H$ by $f_H$. fGN is stationary and from the explicit formula of $f_H$ (cf. Sina{\u\i}~\cite{sin}), we find that $f_H(\lambda)\sim c_H \lambda^{1-2H}$ for $\lambda \downarrow 0$, and that $f_H$ is bounded away from zero elsewhere. Using Lemma~\ref{lemlowerbdstatprocesses},
\[
\lambda_n \bigl(\Cov(\mathbf{Y}_n) \bigr)\geq \biggl(1-
\frac{1}{\pi}\biggr) \inf_{\lambda \in  [1/n,\pi ]} f_H(\lambda)
\gtrsim n^{2H-1}\wedge 1.
\]
From the general formula for the Kullback--Leibler distance between two multivariate normal random variables (or by applying Lemma~\ref{lemKLbd}), we obtain
\begin{eqnarray*}
d_{\KL}(P_1,P_2) &=& \tfrac{1}2 (
\mu_1-\mu_2)^t \Sigma^{-1} (
\mu_1-\mu_2)
\end{eqnarray*}
whenever $P_1$ and $P_2$ denote the probability distributions corresponding to $\mathcal{N}(\mu_1,\Sigma)$ and $\mathcal{N}(\mu_2,\Sigma)$, respectively. This proves the claim.
\end{pf}

As a direct consequence of (57) in \cite{supp}, Lemma~\ref{lemKLbdfGN}, and condition~\textup{(i)} of Theorem~\ref{teomain}, we obtain
the following.

%pr1 #&#
\begin{prop}
\label{propE1E5}
 Given $H\in (0,1)$ suppose that the parameter space $\Theta$ satisfies condition~\textup{(i)} of Theorem~\ref{teomain}. Then,
\begin{eqnarray*}
&&\mE_{1,n}(\Theta)\simeq \mE_{4,n}(\Theta).
\end{eqnarray*}
\end{prop}

%reB.1 #&#
\begin{rem2}
\label{remeasyext}
 The previous proposition can be easily extended to more general stationary noise processes and does not require RKHS theory as only condition~\textup{(i)} of Theorem~\ref{teomain} is involved.
\end{rem2}

%pr2 #&#
\begin{prop}
\label{propE5E2}
Given $H\in (0,1)$ suppose that the parameter space $\Theta$ satisfies condition~\textup{(ii)} of Theorem~\ref{teomain}. Then
\[
\mE_{5,n}(\Theta)\simeq \mE_{2,n}(\Theta).
\]
\end{prop}

\begin{pf}
Recall that $\bbH$ denotes the RKHS associated with $(B_t^H)_{t\in [0,1]}$. From the Moore--Aronszajn theorem, we can conclude $L( \cdot  |   \mathbf{F}_{f,n}) \in \bbH$ since by \eqref{eqLrepresentation} it is a linear combination of functions $K(\cdot, j/n)$.
Condition~\textup{(ii)} of Theorem~\ref{teomain} ensures $F_f\in \bbH$. Define $\bbL_n\subset \bbH$ as the space of functions
\begin{eqnarray*}
&&\sum_{j=1}^n \alpha_j K
\biggl(\cdot, \frac{j}n \biggr)\qquad\mbox{with } (\alpha_1,
\ldots,\alpha_n)^t\in \mathbb{R}^n.
\end{eqnarray*}
From the reproducing property in the Moore--Aronszajn theorem and the interpolation property of $L( \cdot  |   \mathbf{F}_{f,n})$, it follows that $L( \cdot  |   \mathbf{F}_{f,n})\in \bbL_n$ is the projection of $F$ on $\bbL_n$, that is,
\begin{eqnarray*}
&&\langle F, h \rangle_{\bbH} = \bigl\langle L( \cdot |
\mathbf{F}_{f,n}), h \bigr\rangle_{\bbH}\qquad\mbox{for all } h
\in \bbL_n.
\end{eqnarray*}
In particular, $\langle F, L( \cdot  |   \mathbf{F}_{f,n}) \rangle_{\bbH}= \llVert  L( \cdot  |   \mathbf{F}_{f,n})\rrVert  _{\bbH}^2$, and thus
\begin{eqnarray*}
&&\bigl\llVert F - L( \cdot | \mathbf{F}_{f,n}) \bigr\rrVert
_{\bbH}^2 \leq \llVert F - h \rrVert _{\bbH}^2
\qquad\mbox{for all } h \in \bbL_n.
\end{eqnarray*}
Together with Lemma~\ref{lemKLbd}, (57) and condition~\textup{(ii)} in Theorem~\ref{teomain},
\begin{eqnarray*}
\Delta \bigl(\mE_{5,n}(\Theta),\mE_{2,n}(\Theta)
\bigr)^2 &\leq& \sup_{f\in \Theta} d_{\KL}
\bigl(Q_{5,f}^n, Q_f^n \bigr)
\\[-2pt]
&=&\frac{1}2 n^{2-2H} \sup_{f\in \Theta} \bigl\llVert
F - L( \cdot | \mathbf{F}_{f,n}) \bigr\rrVert _{\bbH}^2
\\[-2pt]
&=&\frac{1}2 n^{2-2H} \sup_{f\in \Theta} \inf
_{(\alpha_1,\ldots,\alpha_n)^t\in \mathbb{R}^n} \Biggl\llVert F_f-\sum
_{j=1}^n \alpha_j K \biggl(\cdot,
\frac{j}n \biggr) \Biggr\rrVert _{\bbH}^2
\\[-2pt]
&\rightarrow& 0.
\end{eqnarray*}
This proves the assertion.
\end{pf}

%sC #&#
\section{Proofs for Section~\texorpdfstring{\lowercase{\protect\ref{sec3}}}{3}}\vspace*{-6pt}\label{secproofssecRKHSoffBM}

\begin{pf*}{Proof of Lemma~\ref{lemRiesz}}
From Kadec's $\frac{1}4$-theorem (cf. Young \cite{you}, Theorem~14), we conclude that $(e^{2i\omega_k\cdot})_k$ is a Riesz basis if $\llvert  \omega_k/\pi-k \rrvert  < 1/4$ for all $k\in \mathbb{Z}$. Using
Lemma~\textup{D.1(ii)} and \textup{(iii)} (supplementary material \cite{supp}), we find that $\llvert  \omega_k/\pi-k\rrvert  \leq \frac{1}8\vee \frac{\llvert  1-2H\rrvert  }4<\frac{1}4$ and this proves the claim.
\end{pf*}

%reC.1 #&#
\begin{rem2}
 The constant $\frac{1}4$ in Kadec's $\frac{1}4$-theorem is known to be sharp (cf. \cite{you}, Section~3.3). Since $\omega_k=(k+\frac{1}4(1-2H))\pi+O(1/k)$ by Lemma~\textup{D.1(i)}, the LHS comes arbitrarily close to this upper bound at the boundaries $H\downarrow 0$ and $H\uparrow 1$.
\end{rem2}

\begin{pf*}{Proof of Theorem~\ref{teodzhvzan}}
Recall that $c_H=\sin(\pi H)\Gamma(2H+1)$. In a first step, we prove the identity
%
%eC.1 #&#
\begin{equation}
\frac{c_H}{2\pi}= \frac{2^{4H-3} H\Gamma(H+1/2)\Gamma(3-2H)}{(1-H)\Gamma^2(1-H)\Gamma(3/2-H)}. \label{eqchidentity}
\end{equation}
Application of the replication formula $\Gamma(1-z)\Gamma(z)=\pi/\sin(\pi z)$ for $z=H$ and the duplication formula $\Gamma(2z)=2^{2z-1}\pi^{-1/2}\Gamma(z+1/2)\Gamma(z)$ for $z=H$ and $z=1-H$ gives
\begin{eqnarray*}
\frac{c_H}{2\pi} &=&\frac{\sin(\pi H)\Gamma(2H+1)}{2\pi} =\frac{H\sin(\pi H)\Gamma(2H)}{\pi}
\\
&=&\frac{H\Gamma(2H)}{\Gamma(1-H)\Gamma(H)} =\frac{2^{2H-1} H\Gamma(H+1/2)}{\sqrt{\pi} \Gamma(1-H)}
\\
&=&\frac{2^{2H-2}H\Gamma(H+1/2)\Gamma(3-2H)}{\sqrt{\pi}(1-H)\Gamma(1-H)\Gamma(2-2H)}
=\frac{2^{4H-3} H \Gamma(H+1/2)\Gamma(3-2H)}{(1-H)\Gamma^2(1-H)\Gamma(3/2-H)}.
\end{eqnarray*}
This proves \eqref{eqchidentity}.

Next, let us show that $(\phi_k)_k$ is $L^2(\mu)$-normalized, that is $\llVert  \phi_k\rrVert  _{L^2(\mu)}=1$. This is immediately clear for $k=0$ since (cf. Luke \cite{luk}, Section~13.2)
\[
\int_0^\infty \bigl|  J_{1-H}(\lambda)\bigr|  ^2 \lambda^{-1} \,d\lambda=1/(2-2H).
\]
To compute the normalization constant for $k\neq 0$, the last equality in the proof of \cite{dzh}, Theorem 7.2 gives $\llVert   S_1(2\omega_k, \cdot)\rrVert  _{L^2(\mu)}^{2}= \sigma^{-2}(\omega_k)$, where for $k\neq 0$, using identity \eqref{eqchidentity}, $\sigma^{-2}(\omega_k)=\pi c_H^{-1}2^{2H-2}\llvert  \omega_k\rrvert  ^{2H}J_{-H}^2(\omega_k)$ and $S_1(2\omega_k,2\lambda)=p(H,\omega_k)e^{i(\lambda-\omega_k)}\lambda^HJ_{1-H}(\lambda)/(\lambda-\omega_k)$ with $p(H,\omega_k):= \pi c_H^{-1}2^{2H-2}\omega_k^H\times\break J_{-H}(\omega_k)$. By definition of $\phi_k$ we can write $\phi_k = (\pi/c_H)^{1/2}2^{H-1} \overline{S_1(2\omega_k,\lambda)}/\break p(H,\omega_k)$ and
\begin{eqnarray*}
\llVert \phi_k\rrVert _{L^2(\mu)}^2 &=&
\frac{\pi}{c_H} p\bigl(H,\omega_k\bigr)^{-2}2^{2H-2} \bigl
\llVert S_1(2\omega_k,\cdot) \bigr\rrVert
_{L^2(\mu)}^2
\\
&=&\frac{\pi}{c_H} p(H,\omega_k)^{-2}2^{2H-2}
\sigma^{-2}(\omega_k)
\\
&=&1.
\end{eqnarray*}

Since $\lambda\mapsto\lambda^HJ_{1-H}(\lambda)$ is an odd function, we obtain $S_1(2\omega_{-k},-\lambda)= \overline{S_1(2\omega_k,\lambda)}$ implying  $\phi_k = \psi_{-k}(-\cdot)$ with $\psi_k$ as in Theorem~7.2 of \cite{dzh}. Notice that the space $\mathcal{L}_T$ is defined as the closure of the functions $\overline{\mathcal{F}(\mathbb{I}_t)}$, $t\in [0,1]$, whereas $\bbM$ is the closure of the functions $\mathcal{F}(\mathbb{I}_t)$, $t\in [0,1]$, Therefore, a function $h$ is in $\bbM$ if and only if $h(-\cdot)$ is in $\mathcal{L}_T$. This shows that $\{\phi_k\dvtx  k\in \mathbb{Z}\}$ is a basis of $\bbM$ and that the sampling formula $h= \sum_k a_k h(2\omega_k) \phi_k$ is equivalent to the corresponding result in Theorem 7.2 of \cite{dzh}.

We obtain the expression for $a_0$, using Lebedev \cite{leb}, Formula (5.16.1), $\lim_{\lambda\rightarrow 0} (\lambda/2)^{-\alpha} J_{\alpha}(\lambda) = \Gamma(\alpha+1)^{-1}$, for all $\alpha\geq 0$.

Furthermore, $a_{k}=a_{-k}$ follows from $\omega_{k}=-\omega_{-k}$ and the fact that $a_k^{-1}$ is just a constant times the derivative of $\lambda\mapsto \lambda^H J_{1-H}(\lambda)$ evaluated at $\omega_k$. Since $\lambda\mapsto \lambda^HJ_{1-H}(\lambda)$ is an odd and smooth function, the derivative must be an even function (cf. the remarks after Theorem~\ref{teodzhvzan}) and this gives $a_{k}=a_{-k}$.

To prove \eqref{eqakbehavior}, let us first derive some inequalities. The symbol $\lesssim$ means up to a constant depending on $H$ only.

From the asymptotic expansion of Bessel functions (cf. Gradshteyn and Ryzhik~\cite{gra}, formulas~8.451.1 and 7), $\sum_{r=0}^2 \llvert  J_{r-H}(\lambda)\rrvert  \lesssim  \llvert  \lambda  \rrvert  ^{-1/2}$, for all $\llvert  \lambda\rrvert  \geq \omega_1/2$. Together with Lemma \textup{D.2(ii)} (supplementary material \cite{supp}) applied for $k=0$ and the inequality $\llvert  \mathcal{F}(g_0)\rrvert  \leq \llVert  g_0\rrVert  _{L^1(\mathbb{R})}<\infty$, we find $\llvert  J_{1-H}(\lambda)\rrvert  \lesssim \lambda^{1-H}\wedge \lambda^{-1/2}$. Using Taylor expansion and the recursion formula\break 
$2\frac{d}{d\lambda} J_{1-H}(\lambda)=J_{-H}(\lambda)-J_{2-H}(\lambda)$, for any $k\geq 1$ and any $\lambda \in [\omega_k/2, 2\omega_k ]$,
%
%eC.2 #&#
\begin{equation}
\biggl\llvert \frac{J_{1-H}(\lambda)}{\lambda-\omega_k} \biggr\rrvert \leq \frac{1}2 \sup
_{\xi \in [\omega_k/2,2\omega_k]} \bigl\llvert J_{-H}(\xi)+J_{2-H}(
\xi) \bigr\rrvert \lesssim \llvert \omega_k\rrvert ^{-1/2}.
\label{eqconstevalint3}
\end{equation}

In a second step of the proof, we show that for $k\neq 0$, $G_k(-\infty,\infty)\leq \mathrm{const.} \times \llvert  k\rrvert  ^{1/2-H}$, where
\begin{eqnarray*}
&& G_k(a,b):=\int_a^b \biggl\llvert
\lambda^{-H} \frac{J_{1-H}(\lambda)}{\lambda-\omega_k} \biggr\rrvert \,d\lambda,\qquad k=1,2,
\ldots.
\end{eqnarray*}
Notice that it is enough to prove $G_k(0,\infty)\lesssim \llvert  k\rrvert  ^{1/2-H}$ for $k=1,2,\ldots.$ Decompose $[0,\infty)=[0,\omega_k/2]\cup [\omega_k/2,2\omega_k] \cup [2\omega_k,\infty)$. To bound $G_k(0,\omega_k/2)$, use that $\llvert  \lambda-\omega_k\rrvert  \geq \omega_k/2$ and that $\llvert  J_{1-H}(\lambda)\rrvert  \lesssim \lambda^{1-H}\wedge \lambda^{-1/2}$; to bound $G_k(\omega_k/2,2\omega_k)$, use \eqref{eqconstevalint3}; to bound $G_k(2\omega_k,\infty)$, use that $\llvert  \lambda-\omega_k\rrvert  \geq \lambda/2$ and $\llvert  J_{1-H}(\lambda)\rrvert  \lesssim \lambda^{-1/2}$. Together\vspace*{2pt} with Lemma~\textup{D.1}, this shows that $G_k(0,\infty)\lesssim \llvert  k\rrvert  ^{1/2-H}$.

Next, we show that for $k\neq 0$, $\gamma_{k,H}= e^{i\omega_k}/(\sqrt{2-2H}-1)$,
%
%eC.3 #&#
\begin{equation}
a_k= \gamma_{k,H} a_0+\frac{\sqrt{c_H}\omega_k}{2^H\sqrt{\pi}} \int
e^{i(\omega_k-\lambda)}\frac{\lambda^{H-1}J_{1-H}(\lambda)}{\lambda-\omega_k}\llvert \lambda\rrvert ^{1-2H} \,d
\lambda. \label{eqakidentity}
\end{equation}
Notice that $t^{-1}\mathcal{F}(\mathbb{I}_t)(\lambda)\rightarrow 1$ for $t\rightarrow 0$ and $\lambda$ fixed. Since $\phi_k(\lambda)=\gamma_{k,H}\phi_0(\lambda)+2\omega_k\lambda^{-1}\phi_k(\lambda)$, we have by \eqref{eqQspecexplicitderiv},
\begin{eqnarray*}
a_k &=&\lim_{t\rightarrow 0} \biggl\langle
\phi_k,\frac{1}t \mathcal{F}(\mathbb{I}_t) \biggr
\rangle_{L^2(\mu)}
\\
&=&\gamma_{k,H}a_0+ \lim_{t\rightarrow 0}
\frac{c_H2^{1-2H}}{\pi}\int_{-\infty}^\infty \frac{\omega_k}{\lambda}
\phi_k(2\lambda) \frac{1}t \overline{\mathcal{F}(
\mathbb{I}_t) (2\lambda)} \llvert \lambda\rrvert ^{1-2H} \,d
\lambda.
\end{eqnarray*}
Because of $\llvert  \frac{1}t \overline{\mathcal{F}(\mathbb{I}_t)}\rrvert  \leq 1$, the integrand can be bounded by
\[
\mathrm{const.}\times \bigl\llvert \lambda^{-2H}\phi_k(
\lambda)\bigr\rrvert \lesssim \biggl\llvert \frac{\lambda^{-H}J_{1-H}(\lambda)}{\lambda-\omega_k} \biggr\rrvert.
\]
The $L^1(\mathbb{R})$-norm of this function is smaller than a constant multiple of  $G_k(-\infty,\infty)$ and we may apply dominated convergence, that is, $\lim_{t\rightarrow 0}$ and the integral can be interchanged. The definition of $\phi_k$ gives then \eqref{eqakidentity}.

To prove \eqref{eqakbehavior}, notice that the lower bound follows from \eqref{eqakdef}, \eqref{eqconstevalint3} and Lemma~\textup{D.1} (supplementary material \cite{supp}). For the upper bound, we can restrict ourselves to $k=1,\ldots$ since $a_{k}=a_{-k}$. The statement follows from \eqref{eqakidentity} and $G_k(-\infty,\infty)\lesssim \llvert  k\rrvert  ^{1/2-H}$.
\end{pf*}
\end{appendix}

\section*{Acknowledgement}
The author would like to thank Harrison Zhou for bringing the problem to his attention and helpful discussion. The article was revised based on valuable comments by Yazhen Wang, an
Associate Editor and three anonymous referees.

\begin{supplement}[id=suppA]
\stitle{Asymptotic equivalence for regression under fractional noise\\}
\slink[doi]{10.1214/14-AOS1262SUPP}  %[doi,text={...}] - jei reikia suskaldyti doi
\sdatatype{.pdf}
\sfilename{aos1262\_supp.pdf}
\sdescription{The supplement contains proofs for Section~\ref{secmainresults}, some technical results and a brief summary of the Le Cam distance.}
\end{supplement}

% imsref loaded by linak, 2014-09-10 09:21:35
% imsref loaded by linak, 2014-09-10 10:23:03

\printaddresses
\end{document}